%% file: 2ndConsensus.tex
\documentclass[10pt,twocolumn,twoside,fleqn]{IEEEtran}
\ifCLASSINFOpdf
\else
\fi

\usepackage{graphicx}
\usepackage{amsmath,amsthm}

\usepackage{amssymb}
\usepackage{array}
\usepackage{color}
\usepackage[figuresright]{rotating}
\usepackage{subfigure}
\usepackage{graphics}
\usepackage[hidelinks]{hyperref}
\usepackage[normalem]{ulem}
\usepackage[para,online,flushleft]{threeparttable}
\usepackage{fontenc}
\usepackage[numbers,sort&compress]{natbib}
\usepackage{tabularx}
\usepackage{threeparttable}
\usepackage{makecell}
\usepackage{booktabs}

\allowdisplaybreaks[3]

\newcommand\uu{\boldsymbol{\mathit{u}}}

\newcommand{\Lap}{\ensuremath{\mathbf{L}}}
\newcommand{\A}{\ensuremath{\mathbf{A}}}
\newcommand{\D}{\ensuremath{\mathbf{D}}}

\newtheorem{theorem}{Theorem}[section]
\newtheorem{definition}{Definition}[section]
\newtheorem{lemma}[theorem]{Lemma}

\def\Nsf2{{N}_{\mbox{\scriptsize \rm{2SF}}}}

\begin{document}

\title{Coherence Scaling of Noisy Second-Order Scale-Free Consensus Networks}
\author{Wanyue~Xu, Bin~Wu, Zuobai~Zhang, Zhongzhi~Zhang~\IEEEmembership{Member,~IEEE}, Haibin~Kan, Guanrong~Chen~\IEEEmembership{Fellow, IEEE}
\thanks{The work was supported by the National Natural Science Foundation of China under Grant  61803248, 61872093, U19A2066 and U20B2051, the National Key R \& D Program of China (Nos. 2018YFB1305104 and 2019YFB2101703), Shanghai Municipal Science and Technology Major Project  (No.  2018SHZDZX01), ZJLab, and City University of Hong Kong (Project 7005061). (Corresponding author: Zhongzhi Zhang.)}
\thanks{Wanyue Xu, Bin Wu, Zuobai Zhang, Zhongzhi Zhang, and Haibin Kan are with the Shanghai Key Laboratory of Intelligent Information Processing, School of Computer Science, Fudan University, Shanghai 200433, China; Zhongzhi Zhang and Haibin Kan are also with the Shanghai Engineering Research Institute of Blockchain, Shanghai 200433, China. (e-mail:
	xuwy@fudan.edu.cn;
	binwu11@fudan.edu.cn;
	17300240035@fudan.edu.cn; 
	zhangzz@fudan.edu.cn;
	hbkan@fudan.edu.cn).}
\thanks{Guanrong Chen is with the Department of Electrical Engineering, City University of Hong Kong, Hong Kong SAR, China. (e-mail: eegchen@cityu.edu.hk).}
}
\markboth{IEEE transactions on cybernetics}%
{Xu \MakeLowercase{\textit{et al.}}: Coherence Scaling of Noisy Second-Order Scale-Free Consensus Networks}
%



\IEEEtitleabstractindextext{%
\begin{abstract}
	A striking discovery in the field of network science is that the majority of real networked systems have some universal structural properties. In generally, they are simultaneously sparse, scale-free, small-world, and loopy. In this paper, we investigate the  second-order consensus of dynamic networks with such universal structures subject to white noise at vertices. We focus on the network coherence $H_{\rm SO}$ characterized in terms of the $\mathcal{H}_2$-norm of the vertex systems, which measures the mean deviation of vertex states from their average value. We first study numerically the coherence of some representative real-world networks. We find that their coherence $H_{\rm SO}$ scales sublinearly with the vertex number  $N$. We then study analytically $H_{\rm SO}$ for a class of iteratively growing networks---pseudofractal scale-free webs (PSFWs), and obtain an exact solution to $H_{\rm SO}$, which also increases sublinearly in  $N$, with an exponent much smaller than 1. To  explain the reasons for this sublinear behavior, we finally study  $H_{\rm SO}$ for Sierpin\'ski gaskets, for which $H_{\rm SO}$ grows superlinearly in $N$, with a power exponent much larger than 1.  Sierpin\'ski gaskets have the same number of vertices and edges as the PSFWs, but do not display the scale-free and small-world properties. We thus conclude that the scale-free and small-world, and loopy topologies  are jointly responsible for the observed sublinear scaling of $H_{\rm SO}$.
\end{abstract}

}

\maketitle

\IEEEdisplaynontitleabstractindextext

%
\IEEEpeerreviewmaketitle
\section{Introduction}

\IEEEPARstart{A} {s}   a fundamental problem in interdisciplinary areas ranging from control systems to computer science and physics, consensus has been attracting  extensive attention~\cite{BaDaMa93,SaFaMu07,CaYuReCh12,MoTa14,WuTaCaZh16,ShCaHu17}. It can be applied to various practical scenarios, such as load balancing~\cite{DiFrMo99,AmFrJiVe15}, multi-agent rendezvous~\cite{DiKy07}, UAV flocking~\cite{Sa06}, and sensor networks~\cite{LiRu06,YuChWaYa09,ZhChLiYaGu13}. For multi-agent systems, consensus means that the agents reach an agreement on certain quantities or values, such as load, direction, and pace. However, when the system operates in uncertain environments with noisy disturbances imposed on agents, the system will never reach consensus, with the state of each agent fluctuating around their average. In this case, what we are concerned with is the performance of the system in resisting noise.

The essence of various dynamical networks is the interaction among elements, which can be described by the powerful analytic tool---graphs, where vertices represent the elements and edges represent their relationships~\cite{Ne10}. With this network representation, the interactions of vertex systems are organized into a complex topological structure as a network, characterized by various measurements including degree distribution, average shortest distance, and distribution of cycles or loops of different lengths. These structural properties have striking consequences on the behaviors and performance of dynamical processes running on the networked systems~\cite{Ne03}. In the scenario of networks of agents, the consensus problem has been intensively studied, establishing nontrivial effects of network topological properties on various aspects of the problem, for example, convergence rate~\cite{Ol05,OlTs09,AyBa10, QiZhYiLi19} and robustness to time delay~\cite{XiWa08,MuPaAl10, QiZhYiLi19}, which are determined jointly by the second smallest eigenvalue and the largest  eigenvalue of the Laplacian matrix associated with the graph. 

In addition to convergence rate and time delay, many other interesting quantities about consensus dynamics are also governed by the eigenvalues of the graph Laplacian matrix $\Lap$. For example, for first- and second-order noisy networks without leaders, their network coherence defined in terms  of the system's $\mathcal{H}_2$-norm (i.e., the average of deviations of agent states from the current average value)  is determined by the sum of the reciprocal of the square of each nonzero eigenvalue of $\Lap$~\cite{XiBoKi07,PaBa10,BaJoMiPa12,PaBa14,YiYaZhPa18}. For the first-order consensus problem, the network coherence has been studied for graphs with different structures, including paths~\cite{YoScLe10}, stars~\cite{YoScLe10}, cycles~\cite{YoScLe10}, Vicsek fractal trees and T-fractal trees~\cite{PaBa11,PaBa14}, tori and lattices~\cite{BaJoMiPa12}, Farey graphs~\cite{ZhCo11,YiZhLiCh15}, Koch networks~\cite{YiZhShCh17}, hierarchical graphs~\cite{QiZhYiLi19}, Sierpin\'ski graphs~\cite{ QiZhYiLi19}, and some real-world  networks~\cite{YiZhPa20}. These works revealed nontrivial impacts of the network topology on the behavior of first-order network coherence.

Compared with the first-order setting, the network coherence for second-order consensus problem is relatively rarely studied, in spite of the fact that it can well describe  many practical applications, for example, formation control~\cite{ReAt05} and clock synchronization~\cite{CaZa14}. It has been analyzed only for several special graphs, such as  tori~\cite{BaJoMiPa12}, classic fractals~\cite{PaBa11,PaBa14,QiZhYiLi19}, Koch networks~\cite{YiZhShCh17}, and hierarchical graphs~\cite{QiZhYiLi19}. However, these networks cannot well mimic most real networked  systems, which exhibit universal topological properties~\cite{Ne03}: power-law degree distribution~\cite{BaAl99}, small-world behavior~\cite{WaSt98}, and pattern with cycles at various scales~\cite{RoKiBoBe05,KlSt06}, where a cycle is a path  plus an edge connecting its two ending nodes. It has been shown that the aggregation of these properties has a critical effect on the first-order network coherence~\cite{YiZhPa20}. To date, their effects on second-order network coherence are still largely unknown, which are expected to be quite different from those for the first-order case, since the intrinsic mechanisms governing their dynamics differ significantly.

To fill this gap, in this paper, we study the second-order coherence of noisy consensus on networks with the aforementioned universal properties observed in many real-life networks. 
The main work and contribution of this paper are summarized as follows.

First, we consider the coherence of  scale-free small-world networks with cycles at distinct scales. We show that, for these networks, the second-order coherence scales sublinearly with the number of nodes, which is in sharp contrast to their corresponding first-order coherence that converges to a constant independent of the network size~\cite{YiZhPa20}.

Then, we address the second-order coherence of a family of deterministically iterative networks, called pseudofractal scale-free webs (PSFWs)~\cite{ShLiZh17,ShLiZh18,XiZhCo16}, which display some structural properties similar to those of the real networks studied. By exploiting the self-similarity of the graphs, we establish some recursion relations for the characteristic polynomials of the Laplacian matrices of the PSFWs and their subgraphs at consecutive iterations, based on which we further find exact expressions for the second-order coherence and its leading scaling, which also behaves sublinearly with the network size.

Finally, we show that the sublinear scaling for second-order coherence observed for both real and model networks lies in the composition of scale-free behavior, small-world effect, and the cycles of various scales in the considered networks. For this purpose, we study the second-order coherence of Sierpin\'ski gaskets~\cite{ShLiZh17,ShLiZh18}, which have the same numbers of nodes and edges as the PSFWs but are homogeneous and large-world, an architecture quite different from that of PSFWs. We obtain explicit formulas for the second-order coherence and its dominating behavior, which grows superlinearly with the number of nodes. 

Our results presented in this paper provide insights to understanding the noisy second-order consensus dynamics and have far-reaching implications for the structural design of communication networks.

\section{Preliminaries}\label{coherence.sec}

In this section, we introduce some basic concepts about a graph, its  Laplacian matrix, related distances associated with the eigenvalues and the eigenvectors of the Laplacian matrix, as well as the first-order and second-order noisy consensus problems to be studied. 

\subsection{Graph, Laplacian Matrix and Related Distances}

We use $G=(V, E) $ to denote an undirected connected graph with $N= |V|$ vertices and $M=|E|$ edges, where $V$ is the node set, $E$ is the edge set, and $ |\cdot|$ denotes the cardinality of a set.

The adjacency matrix $\A$ of a  graph $G$ is an $N \times N$ symmetric matrix, representing the adjacent relations of its vertices.  The entry $a_{ij}$ of  $\A$ at row $i$ and column $j$ is defined as follows:  $a_{ij}=a_{ji}=1$ if the vertex pair $(i,j) \in E$, and $a_{ij}=0$ otherwise. Let $\Gamma_i$ be the set of the neighbors for  vertex $i$. 
Then, the degree of vertex $i$ in graph $G$ is defined as  $d_i=\sum_{j=1}^{N} a_{ij}=\sum_{j \in \Gamma_i} a_{ij}$. The average degree of $G$ is $\bar{d}=\frac{1}{N}\sum_{i=1}^{N}d_i=2M/N $. If $\bar{d}$ is a constant, independent of $N$, we call $G$ a sparse graph. For a graph  $G$, its degree matrix $\D$ is a diagonal matrix, with the $i$th diagonal entry equal to $d_i,i=1,2,...,N$.

Let $P(d) $ be the degree distribution of graph  $G$. If   $P(d)\sim d^{-\gamma}$, we call $G$ a scale-free network~\cite{BaAl99}. In a scale-free network, there exist some large-degree nodes, with the maximum-degree vertices called hub vertices, each having degree   $d_{\rm{max}}=N^{1/(\gamma-1)}$. It has been shown~\cite{Ne03} that many real networks are scale-free.

Another important matrix related to a graph $G$ is the Laplacian matrix $\Lap$ defined by $\Lap=\D-\A$~\cite{Me98}. It is an $N \times N$ positive semi-definite matrix with a unique zero eigenvalue and $N-1$ positive  eigenvalues if the graph is connected.  Let $\lambda_i$, $i=1,2,\ldots, N$, be the $N$ eigenvalues of  $\Lap$ rearranged in  ascending order, namely,  $0=\lambda_1<\lambda_2\leq \cdots \leq \lambda_N$. Let $\uu_k$, $ k={1,2,\ldots,N}$, denote the corresponding mutually orthogonal unit eigenvectors, with  the $x$th component being $u_{kx}$, $x={1,2,\ldots,N}$. Using eigenvalues  $\lambda_k$ and their corresponding eigenvectors $\uu_k$,  $ k={1,2,\ldots,N}$, one can define various distances for a graph, such as resistance distance and  biharmonic distance. The resistance distance $\Omega_{ij}$ between two nodes $i$ and $j$ is~\cite{KlRa93}
\begin{equation}\label{rdis}
\Omega_{ij}= \sum_{k=2}^{N}\frac{1}{\lambda_k}(u_{ki}-u_{kj})^2,
\end{equation}
while the biharmonic distance $\Theta_{ij}$ between  $i$ and $j$  is defined as~\cite{LiRuFu10}
\begin{equation}\label{bdis}
\Theta_{ij}=\sum^{N}_{k=2}\frac{1}{\lambda^{2}_{k}}(u_{ki}-u_{kj})^{2} \,.
\end{equation}

The sum of resistance distances over all the $N(N-1)/2$ pairs of vertices in graph $G$ is called its  Kirchhoff index~\cite{KlRa93}, denoted by $R(G)$, which can be represented by all positive eigenvalues of the Laplacian matrix $\Lap$ as~\cite{GhBoSa08}
\begin{equation}
R(G)=\sum\limits_{i,j\in V \atop{i<j}}\Omega_{ij}=\frac{1}{2}\sum_{i,j\in V}\Omega_{ij}=N\,\sum_{i=2}^{N} \frac{1}{\lambda_{i}}\,.\nonumber
\end{equation}
The sum of biharmonic distances over all the $N(N-1)/2$ pairs of vertices in graph $G$ is called its biharmonic index~\cite{YiYaZhPa18}, denoted by $B(G)$. Similarly to $R(G)$,  $B(G)$
can be represented in terms of the $N-1$ non-zero eigenvalues of  $\Lap$: 
\begin{equation}
B(G)=\sum\limits_{{i,j \in V} \atop {i<j}}\Theta_{ij}=N \, \sum\limits^{N}_{i=2}\frac{1}{\lambda^2_i}\,.\nonumber
\end{equation}
It was shown~\cite{TyCoJa18} that $B(G)$ can be expressed in terms of $R(G)$ and the resistance distances of some vertex pairs.

\subsection{Noisy First-Order Consensus Dynamics}

A graph $G$ can be considered  as a multiagent system, where a vertex corresponds to an agent and an edge is associated with  available information flow between two agents. In the first-order consensus network, each agent has a single state. We express the states of the system at time $t$ by an $N$-dimensional real vector $x(t)\in \mathbb{R}^N$, where the $i$th element $x_i(t)$ represents the state of vertex $i$. Every agent adjusts its state according to its local information. In the presence of noise, each agent is subject to stochastic disturbances. For simplicity,  we suppose that every agent is independently influenced by Gaussian  white noise with identical intensity,  which is a stationary and ergodic random process with zero mean and the following fundamental property:  two values of  noise at any pair of times are statistically independent.Then, the evolution of  the system state  can be written in a matrix form as
\begin{equation}\label{HFOdef.eq}
\dot{x} (t)=-\Lap x (t)+w(t),
\end{equation}
where $w(t)=(w_1(t),w_2(t),...,w_N(t))\in \mathbb{R}^N $ is a Gaussian signal with zero-mean and unit variance.


Due to the impact of  noise, the agents will never reach agreement and their states fluctuate around the average value of the current states for all agents. The variance of these fluctuations can be captured by
network coherence, characterized by  the $\mathcal{H}_2$-norm of the system~\cite{YoScLe10}. Without loss
of generality,  we assume that the initial condition $\frac{1}{N}\sum^{N}_{i=1}x_i(0)=0$. The concept of network coherence  represents the extent of the fluctuations~\cite{PaBa11,BaJoMiPa12,PaBa14},

\begin{definition}
	For a graph $G$, the \emph{first-order  network coherence} $H_{\rm FO}$ is defined as the mean steady-state variance of the deviation from the average of the current  agent states:
	\begin{equation}\label{cohDefIn}
	H_{\rm FO} := \frac{1}{N} \lim_{t \rightarrow \infty}    \sum_{i=1}^{N}\mathbf{var}\left \{x_{i}(t) - \frac{1}{N}\sum_{j=1}^{N} x_{j}(t)\right \} \,.
	\end{equation}
\end{definition}

It was shown~\cite{XiBoKi07,BaJoMiPa08,YoScLe10,PaBa11,BaJoMiPa12} that $H_{\rm FO}$ is purely determined by the Kirchhoff index  or the $N-1$ nonzero eigenvalues of $\Lap$, as given by
\begin{equation}\label{Hfo}
H_{\rm FO}=\frac{1}{2N} \sum^{N}_{i=2}\frac{1}{\lambda_i}=\frac{R(G)}{2N^2}.
\end{equation}
$H_{\rm FO}$ measures the performance of the system  robustness to noise. Low $H_{\rm FO}$ corresponds to good robustness, indicating that every agent keeps close to the average of the  current states.

\subsection{Noisy Second-Order Consensus Dynamics}

In the second-order consensus problem, at time $t$, each agent $i$ has two scalar-valued states: $x_{1,i}(t) $ and $x_{2,i}(t)$. Thus, the states of  all agents can be represented by two vectors: position vector $x_1(t) $ and velocity vector $x_2(t)$,  where $x_2(t) $ is the first-order derivative of $x_1(t) $ with respect to time $t$. Different from the first-order case, every agent in the second-order consensus dynamics updates its state by changing the value of $\dot{x_2}(t)$, on the basis of its own states and the states of its neighbors. The noisy second-order consensus system can be described by:
\begin{equation}
\left[ \begin{array}{c}
\dot{x_1}(t)\\
\dot{x_2}(t)
\end{array} \right]=\left[ \begin{array}{cc}
\mathbf{0} & \mathbf{I}\\
-\Lap & -\Lap
\end{array} \right]\left[ \begin{array}{c}
x_1(t)\\
x_2(t)
\end{array} \right]+\left[ \begin{array}{c}
\mathbf{0}\\
\mathbf{I}
\end{array} \right]w(t),
\end{equation}
where vector $w(t) \in \mathbb{R}^{N}$ represents the uncorrelated Gaussian white noise process with zero-mean and unit variance; $\mathbf{0}$ and $\mathbf{I}$ are the $N \times N$ zero matrix and identity matrix, respectively.
 We note that only the variable $x_2(t) $ is subject to disturbances.

The network coherence of the above second-order dynamics only reflects the derivation of state $x_1(t) $ from the average value of the current states of all agents.
\begin{definition}
	For a graph $G$, the \emph{second-order  network coherence} $H_{\rm SO}$ is defined as the mean steady-state variance of the deviation of state $x_1(t) $ from the current average:
	\begin{equation}
	H_{\rm SO}:=\frac{1}{N}\lim_{t\to\infty}\sum\limits^N_{i=1}\mathbf{var}
	\left\{x_{1,i}(t)-\frac{1}{N}\sum\limits^N_{j=1}x_{1,j}(t)  \right\}.
	\end{equation}
\end{definition}
Similarly to $H_{\rm FO}$, $H_{\rm SO}$ is completely determined by the nonzero eigenvalues of the Laplacian matrix~\cite{BaJoMiPa12}.  Specifically, $H_{\rm SO}$ is determined by the  biharmonic index  of the network~\cite{YiYaZhPa18}:
\begin{equation}\label{Hso}
H_{\rm SO}=\frac{1}{2N}\sum\limits^{N}_{i=2}\frac{1}{\lambda_i^2}=\frac{B(G)}{2N^2}.
\end{equation}
A low $H_{\rm SO}$ means that the network structure is robust to random disturbances to the second-order
consensus system.

\subsection{Related Work}

The notion of network coherence was introduced by Bamieh~\emph{et al.} for both noisy first- and second- order consensus dynamics~\cite{BaJoMiPa12}. There are many works focusing on the first-order network coherence. Young~\emph{et al.}~\cite{YoScLe10} derived analytical formulas for first-order network coherence of cycles, paths, and star graphs. Patterson and Bamieh gave exact expressions for the first-order network coherence of some fractal trees~\cite{PaBa14}, as well as tori and lattices~\cite{BaJoMiPa12} of different fractal dimensions. Some co-authors of the present paper presented explicit solutions to the first-order network coherence in Farey graphs~\cite{YiZhLiCh15}, Koch graphs~\cite{YiZhShCh17}, self-similar hierarchical graphs~\cite{QiZhYiLi19} and Sierpi\'nski graphs~\cite{QiZhYiLi19}. These works unveiled some non-trivial effects of network architecture on first-order network coherence. In a recent paper~\cite{YiZhPa20}, the upper and lower bounds of the first-order network were provided for an arbitrary graph, where the lower bound can be approximately reached in most real-world networks.

Relatively to the first-order case, related works about second-order network coherence $H_{\rm SO}$ are much less, with the exception of few  particular graphs, such as tori~\cite{BaJoMiPa12}, fractal trees~\cite{PaBa14}, Koch networks~\cite{YiZhShCh17}, hierarchical graphs~\cite{QiZhYiLi19} and Sierpi\'nski graphs~\cite{QiZhYiLi19}. Very recently, Yi, Zhang, and Patterson~\cite{YiYaZhPa18} established a connection between the biharmonic distance of a graph and its second-order network coherence. They provided exact solutions to second-order network coherence of complete graphs, star graphs, cycles, and paths. However, these studied graphs can not well mimic real-world networks, most of which are sparse, displaying simultaneously scale-free, small-world and loopy structures. Thus far, it has  not been unexplored how the second-order coherence behaves in networks with these general properties. Particularly, there is no exact result about second-order coherence in scale-free, small-world and loopy networks.

In the sequel, we study the second-order network coherence $H_{\rm SO}$ for scale-free, small-world and loopy networks. First, we experimentally study various real-world networks with scale-free small-world structure and loops of different lengths. Then, we derive an exact expression for $H_{\rm SO}$ of a family of scale-free small-world and loopy networks, the PSFWs to be detailed below. We show that their $H_{\rm SO}$ behaves sublinearly with the vertex number $N$. Finally, we obtain an explicit expression for $H_{\rm SO}$ of  loopy Sierpi\'{n}ski gaskets, which are neither scale-free nor small-world, but have the same number of vertices and edges as those of PSFWs. We found that the $H_{\rm SO}$ of  Sierpi\'{n}ski gaskets scales superlinearly with $N$. We argue that the observed sublinear scaling lies in the aggregation of scale-free, small-world, and loopy properties of the studied networks.

\section{Coherence of Some Real  Networks}

In this section, we evaluate the second-order  coherence for 26  real-world networks chosen from the Koblenz Network Collection~\cite{Ku13}, which are  scale-free, small-world, and loopy.  All these realistic networks are typical and representative, including different types of networks such as information networks,  social networks, metabolic
networks, and technological networks.  Table~\ref{tab:network_size} summaries the information of the 26 networks, listed in increasing order of their number of vertices.


\begin{table*}[htbp]
	\normalsize
	\tabcolsep=17pt
	\fontsize{7}{8}\selectfont
	\caption{\label{tab:network_size} Statistics of  26 realistic networks. For a network with $N$ vertices and $M$ edges, we represent the number of vertices and edges in its largest connected component by $N'$ and $M'$, respectively. $\bar{d}'$ represents the average degree of the largest connected component,  equalling $2M'/N'$.  $\gamma$ denotes the power-law exponent. $\bar{l}$ is the average shortest path distance.}
	\resizebox{\textwidth}{!}{
	\begin{tabular}{cccccccc}
		\toprule
		Network & $N$ & $M$ & $N'$ & $M'$ & $\bar{d}'$ & $\gamma$ & $\bar{l}$ \cr
		\midrule
		Karate& 34 &78 & 34 &78 &4.588 & 2.161 &2.408 \cr
		Windsurfers & 43 & 336 & 43 &336 &15.628 & 4.001 &1.671 \cr
		Dolphins &62 & 159 & 62 & 159 &5.129 & 5.001 &3.357 \cr
		Lesmis & 77 & 254 & 77 & 254 &6.597 & 1.521 &2.641 \cr
		Adjnoun & 112 & 425 & 112 & 425 &7.589 &3.621&2.536 \cr
		\scriptsize{ElectronicCircuit(S208)}   & 122 & 189 & 122 & 189 &3.098 & 4.161 &4.928 \cr
		\scriptsize{ElectronicCircuit(S420)}   & 252 & 399 & 252 & 399 &3.167 & 4.021 &5.806 \cr
		Celegansneural & 297 & 2,148 & 297 & 2148 &14.465 & 2.101 &2.455 \cr
		HamsterFull	& 2,426 & 16,631 & 2,000 & 16,098 &16.098 & 2.421 &3.588 \cr
		\scriptsize{WordAdjacency(JAP)}	& 2,704 & 7,998 & 2,698 & 7,995 &5.927 & 2.101 &3.077 \cr
		FacebookNIPS & 4,039 & 88,234 & 4,039 & 88,234 &43.691 & 2.501 &3.693\cr
		GrQc & 5,242 & 14,484 & 4,158 & 13,422 &6.456 & 2.121 &6.049 \cr
		Reactome & 6,327 & 146,160 & 5,973 & 145,778 &48.812 & 1.741 &4.214 \cr
		RouteViews & 6,474 & 12,572  & 6,474 & 12,572 &3.884 & 2.141 &3.705 \cr
		\scriptsize{WordAdjacency(SPA)}	& 7,381 & 44,207 & 7,377 & 44,205 &11.985 & 2.201 &2.778 \cr
		HighEnergy & 7,610 & 15,751 & 5,835 & 13,815 &4.735 &3.441 &7.026 \cr
		HepTh & 9,875 & 25,973 & 8,638 & 24,806 &5.743 & 5.481 &5.945 \cr
		Blogcatalog & 10,312 & 333,983 & 10,312 & 333,983 &64.776 & 2.081 &2.382 \cr
		\scriptsize{PrettyGoodPrivacy} & 10,680 & 24,316 & 10,680 & 24,316 &4.554 & 4.261 &7.486 \cr
		HepPh & 12,006 & 118,489 & 11,204 & 117,619 &20.996 & 2.081 &4.673 \cr
		AstroPh & 18,772 & 198,050 & 17,903 & 196,972 &22.004 & 2.861 &4.194 \cr
		Internet & 22,963 & 48,436 & 22,963 & 48,436 &4.219 & 2.081 &3.842 \cr
		CAIDA & 26,475 & 53,381 & 26,475 & 53,381 &4.033 & 2.101 &3.876\cr
		EnronEmail & 36,692 & 183,831 & 33,696 & 180,811 &10.732 & 1.981 &4.025\cr
		CondensedMatter & 39,577 & 175,692 & 36,458 & 171,735 &9.421 & 3.681 &5.499\cr
		Brightkite & 58,228 & 214,078 & 56,739 & 212,945 &7.353 & 2.501 &4.917 \cr
		\bottomrule
	\end{tabular}}
\end{table*}

\begin{figure}
	\centering
	\includegraphics[width=0.47\textwidth,trim=50 0 50 0]{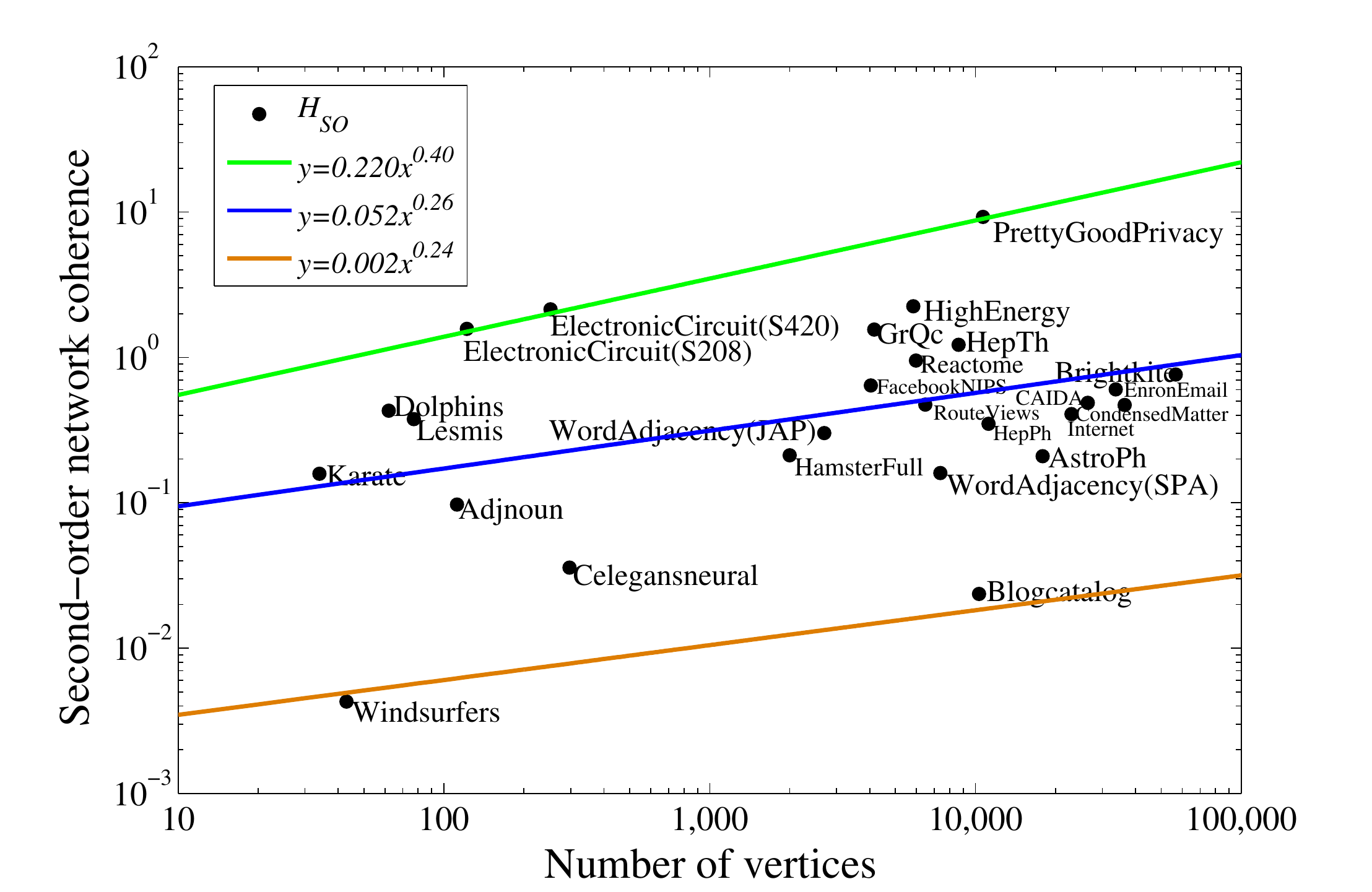}
	\caption{Second-order network coherence  $H_{\rm SO}$ versus  vertex number  in 26 realistic networks on a log–log scale. The solid lines serve as guides.}
	\label{Fig.R}
\end{figure}

Using formula~\eqref{Hso}, we determine the  second-order coherence  $H_{\rm SO}$ for the largest connected component of each studied network, as shown in  Fig.~\ref{Fig.R}. From this figure, we can observe that, for all networks of different sizes, their second-order coherence $H_{\rm SO}$ is approximately  a sublinear  function of their vertex number $N'$, that is, $H_{\rm{SO}} \sim (N')^{\alpha}$, with $0<\alpha<1$. This is in sharp contrast to the first-order coherence  $H_{\rm FO}$, which tends to small constants much less than 1, and is   independent of $N'$~\cite{YiZhPa20}.  

\section{Network Coherence in Pseudofractal Scale-Free Networks}

In this section, we study analytically the coherence for a family of deterministic scale-free model networks, called pseudofractal scale-free webs~\cite{ShLiZh17,ShLiZh18,XiZhCo16}, which display some remarkable properties as observed in many real networks. It is thus expected that the behavior of the network coherence  is similar to that for real networks.


\subsection{Network Construction and Properties}

\begin{figure}
	\centering
	\includegraphics[width=0.30\textwidth]{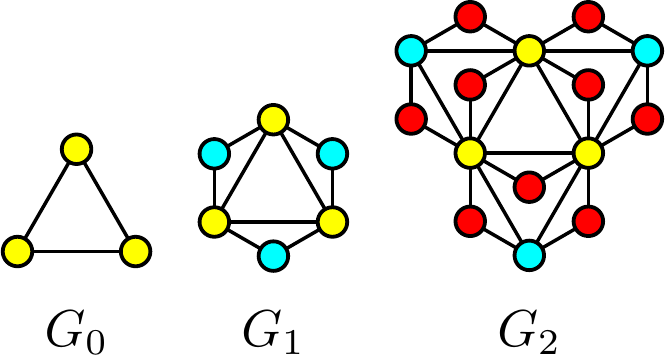}
	\caption{The first three generations of  pseudofractal scale-free webs.}
	\label{Fig.1}
\end{figure}

The pseudofractal scale-free webs (PSFWs) are  generated in an iterative way. We denote by $G_n$, the pseudofractal scale-free network after $n  (n\ge 0)$ iterations. For $n=0$, $G_0$ is a triangle consisting of three vertices and three edges.  When $n\ge 1$, $G_n$ is obtained from $G_{n-1}$ as follows.  Every existing edge in $G_{n-1}$ introduces a new vertex connected to both ends of the edge. Figure~\ref{Fig.1} illustrates the construction process for the first three generations.

Let $N_n$ and  $E_n$ denote, respectively, the number of vertices and the number of edges in $G_n$. It it easy to verify that $N_n=\frac{3^{n+1}+3}{2}$ and $E_n=3^{n+1}$. In  network $G_n$, the three vertices generated at $n=0$ have the largest degree $2^{n+1}$, which are called hub vertices  and  are denoted by $A_n$, $B_n$, and $C_n$, respectively.

\begin{figure}
	\centering
	\includegraphics[width=0.4\textwidth]{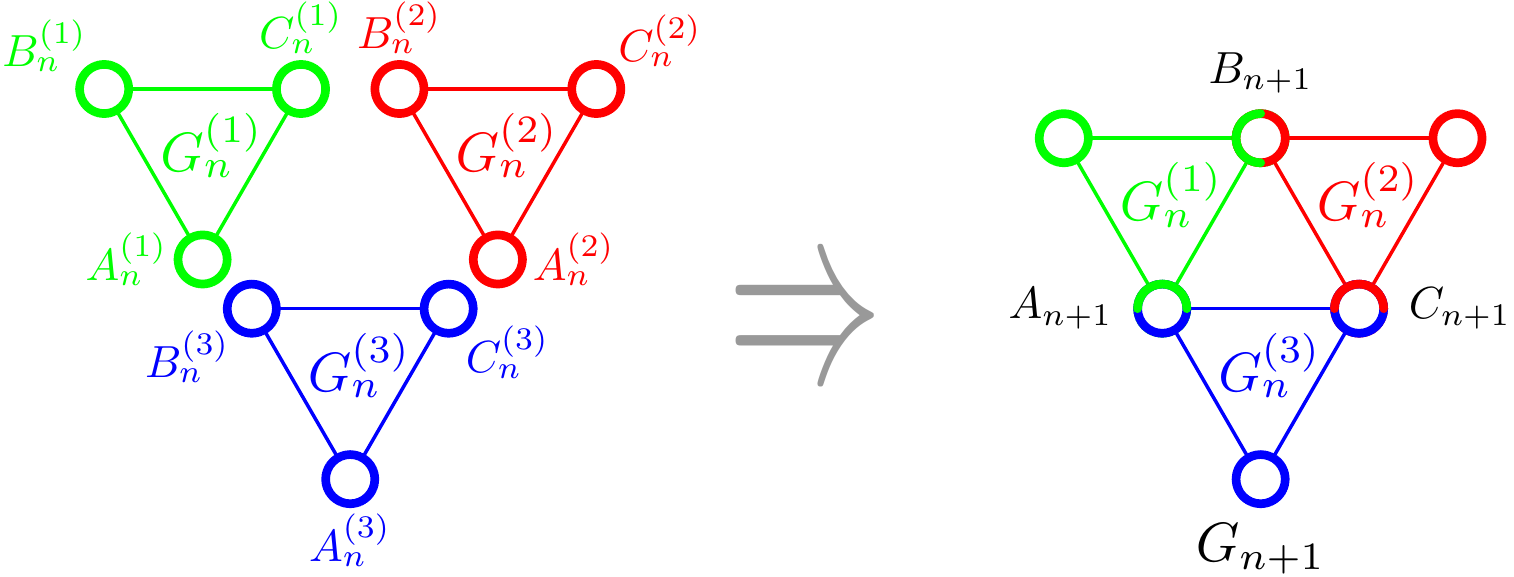}
	\caption{Second construction of pseudofractal scale-free webs, highlighting the self-similar property.}
	\label{Fig.2}
\end{figure}

The  PSFWs are self-similar, which suggests another construction method highlighting their  self-similar structure. This approach creating the  networks is as follows.  Given the $n$th generation network $G_n$,  the $(n+1)$th generation $G_{n+1}$ is obtained by joining three copies of $G_n$ at their hubs, see Fig.~\ref{Fig.2}. Let $G_n^{(\theta)}$, $\theta=1,2,3$, represent the three replicas of $G_n$, and let $A_n^{(\theta)}$, $B_n^{(\theta)}$, and $C_n^{(\theta)}$,  $\theta=1,2,3$, represent, the three hub vertices of  $G_n^{(\theta)}$ ,respectively. Then, $G_{n+1}$ can be generated by merging $G_n^{(\theta)}$, $\theta=1,2,3$, with  $A_n^{(1)}$ and $B_n^{(3)}$ being identified as $A_{n+1}$,  $B_n^{(2)}$ and $C_n^{(1)}$  being identified as $B_{n+1}$, and $A_n^{(2)}$ and $C_n^{(3)}$ being identified as $C_{n+1}$.

The PSFWs exhibit some typical properties of realistic networks. They are scale-free, with the degree distribution $P(d)$ obeying a power law form $P(d)\sim d^{1+\ln 3/\ln 2}$~\cite{DoGoMe02}. They are small-world, with the average distance scaling logarithmically  with $N_n$. Moreover, they are highly clustered, with the average clustering coefficient converging to $\frac{4}{5}$. Finally, they have many cycles of different lengths,  the distribution of which is studied in~\cite{RoKiBoBe05}.

\subsection{Exact Solutions and  Scalings for Network Coherence}

Let $\mathbf{L}_n$ denote the Laplacian matrix of network $G_n$, with a unique zero eigenvalue $\lambda_1(n)$ and $N_n-1$  nonzero eigenvalues $\lambda_2(n)$, $\lambda_3(n)$, $\ldots$, $\lambda_{N_n}(n)$.  Let $H_{\rm FO}(n)$  and  $H_{\rm SO}(n)$  represent, respectively, the first-order network coherence and second-order network coherence of  $G_n$. To determine $H_{\rm FO}(n)$  and  $H_{\rm SO}(n)$, we define two quantities $S_n$ and  $T_n$ by  $S_n=\sum_{i=2}^{N_n} \frac{1}{\lambda_i(n)}$ and $T_n=\sum_{i=2}^{N_n} \frac{1}{\lambda_i^2(n)}$. Then, $H_{\rm FO}(n)=S_n/(2N_n)$  and  $H_{\rm SO}(n)=T_n/(2N_n)$. We next find  $S_n$ and  $T_n$.

\subsubsection{Recursive Relations for Related Polynomials and Quantities}


To determine $S_n$ and  $T_n$, we introduce some quantities.  Let $P_{n}(\lambda)$  denote the characteristic polynomial of matrix ${\bf L}_{n}$, i.e.,
\begin{equation}\label{G7}
P_{n}(\lambda)={\rm det}({\bf L}_{n}-\lambda{\bf I}_{n})\,,
\end{equation}
where ${\bf I}_{n}$ is the $N_n \times N_n$ identity matrix. Let ${\bf Q}_{n}$ be an $(N_n-1)
\times (N_n-1)$ submatrix of $({\bf L}_{n}-\lambda{\bf I}_{n})$,
obtained by removing from $({\bf L}_{n}-\lambda{\bf
	I}_{n})$ the row and column corresponding to a hub vertex in $G_n$. Let ${\bf R}_{n}$
represent a submatrix of $({\bf L}_{n}-\lambda{\bf I}_{n})$ with an
order $(N_n-2) \times (N_n-2)$, obtained from $({\bf
	L}_{n}-\lambda{\bf I}_{n})$ by removing from it  two rows and
columns corresponding to two hub vertices  in $G_n$. Let ${\bf X}_{n}$
represent a submatrix of $({\bf L}_{n}-\lambda{\bf I}_{n})$ with an
order $(N_n-1) \times (N_n-1)$, obtained from $({\bf
	L}_{n}-\lambda{\bf I}_{n})$ by removing from it one row corresponding to a hub vertex and one column corresponding to another hub vertex. Moreover, let $Q_{n}(\lambda)$, $R_{n}(\lambda)$, $X_{n}(\lambda)$ denote, respectively, the determinants of ${\bf Q}_{n}$, ${\bf R}_{n}$ and ${\bf X}_{n}$. As will be shown later, $S_n$ and  $T_n$ can be expressed in terms of the coefficients of some related polynomials. 

\begin{lemma}\label{POLY02}
	For any nonnegative integer $n$,
	\begin{align}
	&P_{n+1}(\lambda)=2 Q_n(\lambda)^3+6 P_n(\lambda) Q_n(\lambda) R_n(\lambda)+\nonumber\\
	&\setlength{\parindent}{4.9em}\indent 9 \lambda  Q_n(\lambda)^2 R_n(\lambda)+3 \lambda  P_n(\lambda) R_n(\lambda)^2+\nonumber\\
	&\setlength{\parindent}{4.9em}\indent 6 \lambda ^2 Q_n(\lambda) R_n(\lambda)^2+\lambda ^3 R_n(\lambda)^3+2 X_n(\lambda)^3\,,\label{D1}\\
	&Q_{n+1}(\lambda)=3 Q_n(\lambda)^2 R_n(\lambda)+P_n(\lambda) R_n(\lambda)^2+\nonumber\\
	&\setlength{\parindent}{4.9em}\indent 4 \lambda  Q_n(\lambda) R_n(\lambda)^2+\lambda ^2 R_n(\lambda)^3\label{D2}\,,\\
	&R_{n+1}(\lambda)=2R_n(\lambda)^2Q_n(\lambda)+\lambda  R_n(\lambda)^3\,,\label{D3}\\
	&X_{n+1}(\lambda)=2 Q_n(\lambda) R_n(\lambda) X_n(\lambda)+\lambda  R_n(\lambda)^2 X_n(\lambda)-\nonumber\\
	&\setlength{\parindent}{4.9em}\indent R_n(\lambda) X_n(\lambda)^2\,.\label{D4}
	\end{align}
\end{lemma}
\begin{IEEEproof}
	By definition, $P_{n+1}$ can be expressed as
	\begin{small}
		\begin{equation}\label{R1}
		P_{n+1}(\lambda)=\left|
		\begin{array}{cccccc}
		2^{n+2}-\lambda&-1&-1&s_n&s_n&0\\
		-1 & 2^{n+2}-\lambda&-1&t_n&0&s_n\\
		-1&-1&2^{n+2}-\lambda&0&t_n&t_n\\
		s_n^T&t_n^T&0&{\bf R}_n&0&0\\
		s_n^T&0&t_n^T&0&{\bf R}_n&0\\
		0&s_n^T&t_n^T&0&0&{\bf R}_n
		\end{array}
		\right|,
		\end{equation}
	\end{small}
where $2^{n+2}$ denotes the degree of the hub vertices $A_{n+1}$, $B_{n+1}$, and $C_{n+1}$ in network $G_{n+1}$; $s_n$ ($t_n$) is a vector of order $N_n-2$ with  $2^{n+1}-1$ nonzero entries $-1$ and $N_n-2^{n+1}-1$ zero entries, in which each $-1$ describes the connection between the hub vertex $A_{n+1}$ ($B_{n+1}$) and vertices belonging to $G_n^{(1)}$, $G_n^{(2)}$ or $G_n^{(3)}$; the superscript $T$ of a vector represents its transpose.
	%
	
	In a similar way,  we obtain
	\begin{equation}\label{R2}
	Q_{n+1}(\lambda)=\left|
	\begin{array}{ccccc}
	2^{n+2}-\lambda  & -1 & t_n & 0 & s_n \\
	-1 & 2^{n+2}-\lambda  & 0 & t_n & t_n \\
	t_n^T & 0 & {\bf R}_n & 0 & 0 \\
	0 & t_n^T & 0 & {\bf R}_n & 0 \\
	s_n^T & t_n^T & 0 & 0 & {\bf R}_n
	\end{array}
	\right|,
	\end{equation}
	\begin{equation}\label{R3}
	R_{n+1}(\lambda)=\left|
	\begin{array}{cccc}
	2^{n+2}-\lambda  & 0 & t_n & t_n \\
	0 & {\bf R}_n & 0 & 0 \\
	t_n^T & 0 & {\bf R}_n & 0 \\
	t_n^T & 0 & 0 & {\bf R}_n
	\end{array}
	\right|,
	\end{equation}
	\begin{equation}\label{R4}
	X_{n+1}(\lambda)=\left|
	\begin{array}{ccccc}
	-1 & -1 & t_n & 0 & s_n \\
	-1 & 2^{n+2}-\lambda  & 0 & t_n & t_n \\
	s_n^T & 0 & {\bf R}_n & 0 & 0 \\
	s_n^T & t_n^T & 0 & {\bf R}_n & 0 \\
	0 & t_n^T & 0 & 0 & {\bf R}_n
	\end{array}
	\right|.
	\end{equation}
	In the sequel, we will show how to derive the recursive relations for $P_{n+1}(\lambda)$, $Q_{n+1}(\lambda)$, $R_{n+1}(\lambda)$, and $X_{n+1}(\lambda)$.
	By the Laplace theorem, we have
	\begin{eqnarray}\label{LE01}
	&&P_{n+1}(\lambda)\nonumber \\
	&&=\left|
	\begin{array}{cccccc}
	2^{n+1}-\lambda  & -1 & 0 & s_{n} & 0 & 0 \\
	-1 & 2^{n+2}-\lambda  & -1 & t_{n} & 0 & s_{n} \\
	-1 & -1 & 2^{n+2}-\lambda  & 0 & t_{n} & t_{n} \\
	s_n^T & t_n^T & 0 & {\bf R}_{n} & 0 & 0 \\
	s_n^T & 0 & t_n^T & 0 & {\bf R}_{n} & 0 \\
	0 & s_n^T & t_n^T & 0 & 0 & {\bf R}_{n}
	\end{array}
	\right|\nonumber\\
	&&+\left|
	\begin{array}{cccccc}
	2^{n+1}-\lambda  & 0 & -1 & 0 & s_{n} & 0 \\
	-1 & 2^{n+2}-\lambda  & -1 & t_{n} & 0 & s_{n} \\
	-1 & -1 & 2^{n+2}-\lambda  & 0 & t_{n} & t_{n} \\
	s_n^T & t_n^T & 0 & {\bf R}_{n} & 0 & 0 \\
	s_n^T & 0 & t_n^T & 0 & {\bf R}_{n} & 0 \\
	0 & s_n^T & t_n^T & 0 & 0 & {\bf R}_{n}
	\end{array}
	\right|\nonumber\\
	&&+\lambda \left|
	\begin{array}{ccccc}
	2^{n+2}-\lambda  & -1 & t_{n} & 0 & s_{n} \\
	-1 & 2^{n+2}-\lambda  & 0 & t_{n} & t_{n} \\
	t_n^T & 0 & {\bf R}_{n} & 0 & 0 \\
	0 & t_n^T & 0 & {\bf R}_{n} & 0 \\
	s_n^T & t_n^T & 0 & 0 & {\bf R}_{n}
	\end{array}
	\right|\,.
	\end{eqnarray}
According to the properties of determinants, it is straightforward to  obtain~\eqref{D1} from~\eqref{LE01}  by using the approach in~\cite{ZhWuLi12}. Similarly, we can derive~\eqref{D2},~\eqref{D3}, and~\eqref{D4}.
\end{IEEEproof}


\subsubsection{Analytical Solutions for Intermediary  Quantities}

Having derived the recursive relations for the above four characteristic polynomial $P_{n}(\lambda)$,
$Q_{n}(\lambda)$, $R_{n}(\lambda)$ and $X_{n}(\lambda)$, we now determine the coefficients of $P_n(\lambda)$. Define $p_{n}^{(i)}$ ($0 \leq i \leq 2$) as the coefficient of the term $\lambda^i$ corresponding to $\frac{P_n(\lambda)}{\lambda}$. Then, $p_{n}^{(0)}$ denotes the constant item, $p_{n}^{(1)}$ and $p_{n}^{(2)}$ are the coefficients of the terms with degree 1 and 2, respectively. According to Vieta's formulas, we obtain
\begin{align}
&S_n=\sum_{i=2}^{N_n} \frac{1}{\lambda_i(n)}=-\frac{p_n^{(1)}}{p_n^{(0)}}\,,\label{S01}\\
&T_n=\sum_{i=2}^{N_n} \frac{1}{\lambda_i^2(n)}\nonumber\\
&\quad=\left(\sum_{i=2}^{N_n} \frac{1}{\lambda_i(n)}\right)^2-2\sum_{2\leq i < j \leq N_n}\frac{1}{\lambda_i (n)\lambda_j(n)}\nonumber\\
&\quad=S_n^2-2\frac{p_n^{(2)}}{p_n^{(0)}}\,.\label{S02}
\end{align}
Thus, the problem of determining $S_n$ and $T_n$ is reduced to determining  $p_n^{(0)}$, $p_n^{(1)}$, and $p_n^{(2)}$. In order to find these three coefficients, we introduce some additional quantities.  Let $q_{n}^{(i)}$, $r_{n}^{(i)}$, $x_{n}^{(i)}$, $0 \leq i \leq 3$, be the coefficients of term $\lambda^i$ corresponding to $Q_n(\lambda)$ ,$R_n(\lambda)$, $X_n(\lambda)$, respectively.
\begin{lemma} \label{LM01}
	For any nonnegative integer $n$,
	\begin{align}
	&p_n^{(0)}=-2^{\frac{1}{4} \left(-7+3^{1+n}-2 n\right)} 3^{\frac{1}{4} \left(5+3^{1+n}+2 n\right)} \left(1+3^n\right)\,,\label{A11}\\
	&q_n^{(0)}=2^{-\frac{3}{4}+\frac{3^{1+n}}{4}-\frac{n}{2}} 3^{\frac{1}{4}+\frac{3^{1+n}}{4}+\frac{n}{2}}\,,\label{A12}\\
	&r_n^{(0)}=2^{\frac{1}{4}+\frac{3^{1+n}}{4}+\frac{n}{2}} 3^{-\frac{3}{4}+\frac{3^{1+n}}{4}-\frac{n}{2}}\,,\label{A13}\\
	&x_n^{(0)}=-2^{-\frac{3}{4}+\frac{3^{1+n}}{4}-\frac{n}{2}} 3^{\frac{1}{4}+\frac{3^{1+n}}{4}+\frac{n}{2}}\,,\label{A14}\\
	%
	&p_n^{(1)}=\frac{1}{7} 2^{\frac{1}{4} \left(-15+3^{1+n}-2 n\right)} 3^{\frac{1}{4} \left(1+3^{1+n}-2 n\right)}\times\nonumber\\
	&\setlength{\parindent}{3.2em}\indent \big(25 \times 2^n-7 \times 3^n+8 \times 3^{1+2 n}+\nonumber\\
	&\setlength{\parindent}{3.2em}\indent 25 \times 3^{1+3 n}+5 \times 6^{1+n}-35 \times 18^n\big)\,,\label{A21}\\
	&q_n^{(1)}=\frac{1}{7} 2^{\frac{1}{4} \left(-11+3^{1+n}-2 n\right)} 3^{\frac{1}{4} \left(-3+3^{1+n}-2 n\right)}\times\nonumber\\
	&\setlength{\parindent}{3.2em}\indent  \left(-11 \times 2^{2+n}+7 \times 3^n-25 \times 3^{1+2 n}\right)\,,\label{A22}\\
	&r_n^{(1)}=\frac{1}{7} 2^{\frac{1}{4} \left(-7+3^{1+n}+2 n\right)} 3^{\frac{1}{4} \left(-7+3^{1+n}-6 n\right)}\times\nonumber\\
	&\setlength{\parindent}{3.2em}\indent \left(3 \times 2^{2+n}-25 \times 3^{1+2 n}+7\times 3^n \left(-1+2^{2+n}\right)\right)\,,\label{A23}\\
	&x_n^{(1)}=\frac{1}{7} 2^{\frac{1}{4} \left(-11+3^{1+n}-2 n\right)} 3^{\frac{1}{4} \left(-3+3^{1+n}-2 n\right)} \times\nonumber\\
	&\setlength{\parindent}{3.2em}\indent  \left(2^{1+n}-7 \times 3^n+25 \times 3^{1+2 n}-7 \times 6^{1+n}\right)\,,\label{A24}\\
	%
	&p_n^{(2)}=\frac{1}{5635}2^{\frac{1}{4} \left(-27+3^{1+n}-2 n\right)} 3^{\frac{1}{4} \left(-7+3^{1+n}-6 n\right)}\times\nonumber\\
	&\setlength{\parindent}{3.2em}\indent \big(-41\times 2^{7+2 n} 3^{1+n}-9775\times 2^{1+n} 3^{3+2 n}+\nonumber\\
	&\setlength{\parindent}{3.2em}\indent 129283\times 3^{1+3 n}-71875\times 3^{3+5 n}+9039\times 4^{2+n}-\nonumber\\
	&\setlength{\parindent}{3.2em}\indent 20125\times 6^{1+n}+93541\times 9^n+79373\times 4^{2+n} 9^n+\nonumber\\
	&\setlength{\parindent}{3.2em}\indent 147163\times 9^{1+2 n}+100625\times 2^{1+n} 9^{1+2 n}-\nonumber\\
	&\setlength{\parindent}{3.2em}\indent 64975\times 54^{1+n}\big)\,,\label{A31}\\
	&q_n^{(2)}=\frac{1}{5635}2^{\frac{1}{4} \left(-23+3^{1+n}-2 n\right)} 3^{\frac{1}{4} \left(-11+3^{1+n}-6 n\right)}\times\nonumber\\
	&\setlength{\parindent}{3.2em}\indent \big(8855 \times 2^{3+n} 3^{1+n}-1127 \times 2^{7+2 n} 3^{1+n}+\nonumber\\
	&\setlength{\parindent}{3.2em}\indent 71875 \times 3^{3+4 n}-18819 \times 4^{2+n}-93541 \times 9^n-\nonumber\\
	&\setlength{\parindent}{3.2em}\indent 61985 \times 4^{2+n} 9^n+31625 \times 2^{3+n} 9^{1+n}-\nonumber\\
	&\setlength{\parindent}{3.2em}\indent 36596 \times 27^{1+n}\big)\,,\label{A32}\\
	&r_n^{(2)}=\frac{1}{5635}2^{\frac{1}{4} \left(-19+3^{1+n}+2 n\right)} 3^{\frac{1}{4} \left(-15+3^{1+n}-10 n\right)}\times\nonumber\\
	&\setlength{\parindent}{3.2em}\indent \big(18873\times 2^{3+2 n}+161\times 2^{4+2 n} 3^{3+n}-\nonumber\\
	&\setlength{\parindent}{3.2em}\indent 115\times 2^{9+n} 3^{1+2 n}-29288\times 3^{2+3 n}-\nonumber\\
	&\setlength{\parindent}{3.2em}\indent 20125\times 2^{3+n} 3^{2+3 n}+71875\times 3^{3+4 n}-\nonumber\\
	&\setlength{\parindent}{3.2em}\indent 805\times 6^{3+n}+127351\times 9^n-\nonumber\\
	&\setlength{\parindent}{3.2em}\indent 28175\times 2^{3+2 n} 9^n\big)\,,\label{A33}\\
	&x_n^{(2)}=\frac{1}{5635}2^{\frac{1}{4} \left(-23+3^{1+n}-2 n\right)} 3^{\frac{1}{4} \left(-11+3^{1+n}-6 n\right)}\times\nonumber\\
	&\setlength{\parindent}{3.2em}\indent \big(1413\times 2^{3+2 n}-805\times 2^{2+n} 3^{1+n}+\nonumber\\
	&\setlength{\parindent}{3.2em}\indent 1771\times 2^{4+2 n} 3^{1+n}-71875\times 3^{3+4 n}+\nonumber\\
	&\setlength{\parindent}{3.2em}\indent 93541\times 9^n-28175\times 2^{3+2 n} 9^n-\nonumber\\
	&\setlength{\parindent}{3.2em}\indent 8165\times 2^{4+n} 9^{1+n}+36596\times 27^{1+n}+\nonumber\\
	&\setlength{\parindent}{3.2em}\indent 20125\times 2^{2+n} 27^{1+n}\big)\,.\label{A34}
	\end{align}
\end{lemma}
\begin{IEEEproof}
	From~\eqref{D1}-\eqref{D4}, by using an  approach similar to that in~\cite{ZhWuLi12}, it is not difficult to derive the following recursive relations governing the above-defined coefficients:
	\begin{align}
	&p_{n+1}^{(0)}=6 p_n^{(0)} q_n^{(0)} r_n^{(0)}+9 [q_n^{(0)}]^2 r_n^{(0)}+6 [q_n^{(0)}]^2 q_n^{(1)}+\nonumber\\
	&\setlength{\parindent}{3.4em}\indent 6 [x_n^{(0)}]^2 x_n^{(1)}\,,\label{C11} \\
	&q_{n+1}^{(0)}=3 [q_n^{(0)}]^2 r_n^{(0)}\,,\label{C12}\\
	&r_{n+1}^{(0)}=2 q_n^{(0)} [r_n^{(0)}]^2\,,\label{C13} \\
	&x_{n+1}^{(0)}=2 q_n^{(0)} r_n^{(0)} x_n^{(0)}-r_n^{(0)}[x_n^{(0)}]^2\,,\label{C14}\\
	%
	&p_{n+1}^{(1)}=3 p_n^{(0)} [r_n^{(0)}]^2+6 q_n^{(0)} [r_n^{(0)}]^2+6 q_n^{(0)} r_n^{(0)} p_n^{(1)}+\nonumber\\
	&\setlength{\parindent}{3.4em}\indent 6 p_n^{(0)} r_n^{(0)} q_n^{(1)}+18 q_n^{(0)} r_n^{(0)} q_n^{(1)}+6 q_n^{(0)} [q_n^{(1)}]^2+\nonumber\\
	&\setlength{\parindent}{3.4em}\indent 6 p_n^{(0)} q_n^{(0)} r_n^{(1)}+9 [q_n^{(0)}]^2 r_n^{(1)}+6 x_n^{(0)} [x_n^{(1)}]^2+\nonumber\\
	&\setlength{\parindent}{3.4em}\indent [6 q_n^{(0)}]^2 q_n^{(2)}+6 [x_n^{(0)}]^2 x_n^{(2)}\,,\label{C21}\\
	&q_{n+1}^{(1)}=p_n^{(0)} [r_n^{(0)}]^2+4 q_n^{(0)} [r_n^{(0)}]^2+6 q_n^{(0)} r_n^{(0)} q_n^{(1)}+\nonumber\\
	&\setlength{\parindent}{3.4em}\indent 3 [q_n^{(0)}]^2 r_n^{(1)}\,,\label{C22}\\
	&r_{n+1}^{(1)}=[r_n^{(0)}]^3+2 [r_n^{(0)}]^2 q_n^{(1)}+4 q_n^{(0)} r_n^{(0)} r_n^{(1)}\,,\label{C23}\\
	&x_{n+1}^{(1)}=[r_n^{(0)}]^2 x_n^{(0)}+2 r_n^{(0)} x_n^{(0)} q_n^{(1)}+2 q_n^{(0)} x_n^{(0)} r_n^{(1)}-\nonumber\\
	&\setlength{\parindent}{3.4em}\indent [x_n^{(0)}]^2 r_n^{(1)}+2 q_n^{(0)} r_n^{(0)} x_n^{(1)}-2 r_n^{(0)} x_n^{(0)} x_n^{(1)}\,,\label{C24}\\
	%
	&p_{n+1}^{(2)}=[r_n^{(0)}]^3+3 [r_n^{(0)}]^2 p_n^{(1)}+6 [r_n^{(0)}]^2 q_n^{(1)}+6 r_n^{(0)} p_n^{(1)}\nonumber\\
	&\setlength{\parindent}{3.4em}\indent  q_n^{(1)}+9 r_n^{(0)} [q_n^{(1)}]^2+2 [q_n^{(1)}]^3+6 p_n^{(0)} r_n^{(0)} r_n^{(1)}+\nonumber\\
	&\setlength{\parindent}{3.4em}\indent 12 q_n^{(0)} r_n^{(0)} r_n^{(1)}+6 q_n^{(0)} p_n^{(1)} r_n^{(1)}+6 p_n^{(0)} q_n^{(1)} r_n^{(1)}+\nonumber\\
	&\setlength{\parindent}{3.4em}\indent 18 q_n^{(0)} q_n^{(1)} r_n^{(1)}+2 [x_n^{(1)}]^3+6 q_n^{(0)} r_n^{(0)} p_n^{(2)}+\nonumber\\
	&\setlength{\parindent}{3.4em}\indent 6 p_n^{(0)} r_n^{(0)} q_n^{(2)}+18 q_n^{(0)} r_n^{(0)} q_n^{(2)}+12 q_n^{(0)} q_n^{(1)} q_n^{(2)}+\nonumber\\
	&\setlength{\parindent}{3.4em}\indent 6 p_n^{(0)} q_n^{(0)} r_n^{(2)}+9 [q_n^{(0)}]^2 r_n^{(2)}+12 x_n^{(0)} x_n^{(1)} x_n^{(2)}+\nonumber\\
	&\setlength{\parindent}{3.4em}\indent 6 [q_n^{(0)}]^2 q_n^{(3)}+6 [x_n^{(0)}]^2 x_n^{(3)}\,,\label{C31}\\
	&q_{n+1}^{(2)}=[r_n^{(0)}]^3+[r_n^{(0)}]^2 p_n^{(1)}+4 [r_n^{(0)}]^2 q_n^{(1)}+3 r_n^{(0)}\nonumber\\
	&\setlength{\parindent}{3.4em}\indent  [q_n^{(1)}]^2+2 p_n^{(0)} r_n^{(0)} r_n^{(1)}+8 q_n^{(0)} r_n^{(0)} r_n^{(1)}+6 q_n^{(0)} \nonumber\\
	&\setlength{\parindent}{3.4em}\indent q_n^{(1)} r_n^{(1)}+6 q_n^{(0)} r_n^{(0)} q_n^{(2)}+3 [q_n^{(0)}]^2 r_n^{(2)}\,,\label{C32}\\
	&r_{n+1}^{(2)}=3 [r_n^{(0)}]^2 r_n^{(1)}+4 r_n^{(0)} q_n^{(1)} r_n(1)+2 q_n^{(0)} [r_n^{(1)}]^2+\nonumber\\
	&\setlength{\parindent}{3.4em}\indent 2 [r_n^{(0)}]^2 q_n^{(2)}+4 q_n^{(0)} r_n^{(0)} r_n^{(2)}\,,\label{C33}\\
	&x_{n+1}^{(2)}=2 r_n^{(0)} x_n^{(0)} r_n^{(1)}+2 x_n^{(0)} q_n^{(1)} r_n^{(1)}+[r_n^{(0)}]^2 x_n^{(1)}+\nonumber\\
	&\setlength{\parindent}{3.4em}\indent 2 r_n^{(0)} q_n^{(1)} x_n^{(1)}+2 q_n^{(0)} r_n^{(1)} x_n^{(1)}-2 x_n^{(0)} r_n^{(1)} x_n^{(1)}-\nonumber\\
	&\setlength{\parindent}{3.4em}\indent r_n^{(0)} [x_n^{(1)}]^2+2 r_n^{(0)} x_n^{(0)} q_n^{(2)}+2 q_n^{(0)} x_n^{(0)} r_n^{(2)}-\nonumber\\
	&\setlength{\parindent}{3.4em}\indent [x_n^{(0)}]^2 r_n^{(2)}+2 q_n^{(0)} r_n^{(0)} x_n^{(2)}-2 r_n^{(0)} x_n^{(0)} x_n^{(2)}\,,\label{C34}\\
	%
	&q_{n+1}^{(3)}=3 [r_n^{(0)}]^2 r_n^{(1)}+2 r_n^{(0)} p_n^{(1)} r_n^{(1)}+8 r_n^{(0)} q_n^{(1)} r_n^{(1)}+\nonumber\\
	&\setlength{\parindent}{3.4em}\indent 3 [q_n^{(1)}]^2 r_n^{(1)}+p_n^{(0)} [r_n^{(1)}]^2+4 q_n^{(0)} [r_n^{(1)}]^2+\nonumber\\
	&\setlength{\parindent}{3.4em}\indent [r_n^{(0)}]^2 p_n^{(2)}+4 [r_n^{(0)}]^2 q_n^{(2)}+6 r_n^{(0)} q_n^{(1)} q_n^{(2)}+\nonumber\\
	&\setlength{\parindent}{3.4em}\indent 6 q_n^{(0)} r_n^{(1)} q_n^{(2)}+2 p_n^{(0)} r_n^{(0)} r_n^{(2)}+8 q_n^{(0)} r_n^{(0)} r_n^{(2)}+\nonumber\\
	&\setlength{\parindent}{3.4em}\indent 6 q_n^{(0)} q_n^{(1)} r_n^{(2)}+6 q_n^{(0)} r_n^{(0)} q_n^{(3)}+3 [q_n^{(0)}]^2 r_n^{(3)}\,,\label{C42}\\
	&r_{n+1}^{(3)}=3 r_n^{(0)} [r_n^{(1)}]^2+2 q_n^{(1)} [r_n^{(1)}]^2+4 r_n^{(0)} r_n^{(1)} q_n^{(2)}+\nonumber\\
	&\setlength{\parindent}{3.4em}\indent 3 [r_n^{(0)}]^2 r_n^{(2)}+4 r_n^{(0)} q_n^{(1)} r_n^{(2)}+4 q_n^{(0)} r_n^{(1)} r_n^{(2)}+\nonumber\\
	&\setlength{\parindent}{3.4em}\indent 2 [r_n^{(0)}]^2 q_n^{(3)}+4 q_n^{(0)} r_n^{(0)} r_n^{(3)}\,,\label{C43}\\
	&x_{n+1}^{(3)}=x_n^{(0)} [r_n^{(1)}]^2+2 r_n^{(0)} r_n^{(1)} x_n^{(1)}+2 q_n^{(1)} r_n^{(1)} x_n^{(1)}-\nonumber\\
	&\setlength{\parindent}{3.4em}\indent r_n^{(1)} [x_n^{(1)}]^2+2 x_n^{(0)} r_n^{(1)} q_n^{(2)}+2 r_n^{(0)} x_n^{(1)} q_n^{(2)}+\nonumber\\
	&\setlength{\parindent}{3.4em}\indent 2 r_n^{(0)} x_n^{(0)} r_n^{(2)}+2 x_n^{(0)} q_n^{(1)} r_n^{(2)}+2 q_n^{(0)} x_n^{(1)} r_n^{(2)}-\nonumber\\
	&\setlength{\parindent}{3.4em}\indent 2 x_n^{(0)} x_n^{(1)} r_n^{(2)}+[r_n^{(0)}]^2 x_n^{(2)}+2 r_n^{(0)} q_n^{(1)} x_n^{(2)}+\nonumber\\
	&\setlength{\parindent}{3.4em}\indent 2 q_n^{(0)} r_n^{(1)} x_n^{(2)}-2 x_n^{(0)} r_n^{(1)} x_n^{(2)}-2 r_n^{(0)} x_n^{(1)} x_n^{(2)}+\nonumber\\
	&\setlength{\parindent}{3.4em}\indent 2 r_n^{(0)} x_n^{(0)} q_n^{(3)}+2 q_n^{(0)} x_n^{(0)} r_n^{(3)}-[x_n^{(0)}]^2 r_n^{(3)}+\nonumber\\
	&\setlength{\parindent}{3.4em}\indent 2 q_n^{(0)} r_n^{(0)} x_n^{(3)}-2 r_n^{(0)} x_n^{(0)} x_n^{(3)}\,.\label{C44}
	\end{align}
	With the initial values $p_0^{(0)}$=-9, $q_0^{(0)}$=3, $r_0^{(0)}$=2, $x_0^{(0)}$=3, $p_0^{(1)}$=6, $q_0^{(1)}$=-4, $r_0^{(1)}$=-1, $x_0^{(1)}$=1, $p_0^{(2)}$=-1, $q_0^{(2)}$=-1, $r_0^{(2)}$=0, and $x_0^{(2)}$=0, \eqref{C11}-\eqref{C44} can be solved to obtain \eqref{A11}-\eqref{A34}.
	Due to the space limitation, below  we only derive $q^{(0)}_{n}$, $r^{(0)}_{n}$ and $x^{(0)}_{n}$, the other quantities can be obtained in a similar way. 

	\eqref{C13} can be rewritten as
	\begin{equation}\label{r0r}
	q^{(0)}_{n}=\frac{r^{(0)}_{n+1}}{2\left[r^{(0)}_{n}\right]^2}.
	\end{equation} 
	Inserting \eqref{r0r} into \eqref{C12} leads to
	\begin{equation}\label{r0re}
	\frac{r^{(0)}_{n+2}}{\left[r^{(0)}_{n+1}\right]^3}=\frac{3r^{(0)}_{n+1}}{2\left[r^{(0)}_{n}\right]^3},
	\end{equation}
	which provides an explicit  recursive relation governing $r^{(0)}_{n},\ r^{(0)}_{n+1},$ and $r^{(0)}_{n+2}$.
		
	We now derive a closed-form expression for $r^{(0)}_{n}$. To this end, we introduce an intermediary quantity $k_n$, defined as
		\begin{equation}\label{rtok}
		k_n=\frac{r^{(0)}_{n}}{\left[r^{(0)}_{n-1}\right]^3},
		\end{equation}
		which, together with \eqref{r0re}, yields 
		\begin{equation}\label{kre}
		k_{n+1}=\frac{3}{2}k_{n}.
		\end{equation}
		Using the initial condition $k_1=r^{(0)}_{1}/\left[r^{(0)}_{0}\right]^3=3$, \eqref{kre} is solved to give
		\begin{equation}\label{kso}
		k_n=2\left(\frac{3}{2}\right)^n.
		\end{equation}
		With this exact result of $k_n$,  \eqref{rtok} is recast as
		\begin{equation}\label{lnr}
		\ln{r^{(0)}_{n}}=3\ln{r^{(0)}_{n-1}}+\ln{k^{(0)}_{n}}.
		\end{equation}
		Considering $\ln{r^{(0)}_{0}}=\ln{2}$ and the expression for $k_n$ given in \eqref{kso}, \eqref{lnr} is solved inductively
		to yield
		\begin{equation}
		\ln{r^{(0)}_{n}}=3^n\ln{r^{(0)}_{0}}+\sum_{i=0}^{n-1}3^i\ln{\left[2\left(\frac{3}{2}\right)^{i+1}\right]},
		\end{equation} 
		which implies 
		\begin{equation}
		r_n^{(0)}=2^{\frac{1}{4}+\frac{3^{1+n}}{4}+\frac{n}{2}} 3^{-\frac{3}{4}+\frac{3^{1+n}}{4}-\frac{n}{2}}\,,\nonumber
		\end{equation}
		as shown in~\eqref{A13}. 
		
	After deriving $r^{(0)}_{n}$, we continue to calculate $q^{(0)}_{n}$. Plugging~\eqref{A13} into \eqref{r0r} gives 
		\begin{equation}\label{q0so}
		q_n^{(0)}=2^{-\frac{3}{4}+\frac{3^{1+n}}{4}-\frac{n}{2}} 3^{\frac{1}{4}+\frac{3^{1+n}}{4}+\frac{n}{2}}.\nonumber
		\end{equation}
	In this way, we obtain~\eqref{A12}. Finally, we determine $x^{(0)}_{n}$. By inserting  \eqref{C12} into \eqref{C14}, one obtains 
		\begin{equation}\label{x0re}
		\frac{x^{(0)}_{n+1}}{q^{(0)}_{n+1}}=\frac{2x^{(0)}_{n}}{3q^{(0)}_{n}}-\frac{1}{3}\left[\frac{x^{(0)}_{n}}{q^{(0)}_{n}}\right]^2.
		\end{equation} 
		In order to obtain the exact expression for $x^{(0)}_{n}$, we introduce another  quantity $y_n$ defined by
		\begin{equation}\label{yde}
		y_n=\frac{x^{(0)}_{n}}{q^{(0)}_{n}}.
		\end{equation}
		Then, \eqref{x0re} can be rewritten as
		\begin{equation}\label{yre}
		y_{n+1}=\frac{2}{3}y_{n}-\frac{1}{3}y_{n}^2,
		\end{equation}
		which, under the initial condition $y_0=x^{(0)}_{0}/q^{(0)}_{0}=-1$, is  solved to yield
		\begin{equation}\label{yso}
		y_n=-1.
		\end{equation} 
		Combine \eqref{yso}, \eqref{yde} and $r_n^{(0)}$, we obtain
		\begin{equation}
		x_n^{(0)}=-2^{-\frac{3}{4}+\frac{3^{1+n}}{4}-\frac{n}{2}} 3^{\frac{1}{4}+\frac{3^{1+n}}{4}+\frac{n}{2}},\nonumber
		\end{equation}
		as provided by~\eqref{A14}.
\end{IEEEproof}

\subsubsection{Explicit Expression and Behavior  for Network Coherence}

With the above-obtained quantities, we can get  accurate solutions for both the first-order and the second-order network coherence of  network $G_n$, from which we can further reveal their asymptotical behaviors.

\begin{theorem}\label{TH01}
	For the PSFW $G_n$ with $n \geq 1$, the first-order coherence $H_{\rm FO}(n)$ and the second-order coherence $H_{\rm SO}(n)$ of the system with dynamics in \eqref{Hfo} and \eqref{Hso} are
	\begin{small}
		\begin{align}
		&H_{\rm FO}(n)=\frac{1}{28 \left(1+3^n\right)^2 3^{2+n}}\times\big(25\times 2^n-\nonumber\\
		&\setlength{\parindent}{4.8em}\indent 7\times 3^n+8\times 3^{1+2 n}+25\times 3^{1+3 n}+\nonumber\\
		&\setlength{\parindent}{4.8em}\indent 5\times 6^{1+n}-35\times 18^n\big)\,,\label{T1}\\
		&H_{\rm SO}(n)=\frac{3^{-4-2 n}}{90160 \left(1+3^n\right)^3}\times \big(69538\times 3^{2+5 n}+360249\times\nonumber\\
		&\setlength{\parindent}{4.8em}\indent 4^n+ 35\times 2^{2+n} 3^{1+n} \left(-575+1539\times 2^n\right)+\nonumber\\
		&\setlength{\parindent}{4.8em}\indent 322\times 27^n \left(1135-3225\times 2^{1+n}+847\times 2^{1+2 n}\right)+\nonumber\\
		&\setlength{\parindent}{4.8em}\indent 3^{1+4 n} \left(516262-60375\times 2^{2+n}+140875\times 4^n\right)+\nonumber\\
		&\setlength{\parindent}{4.8em}\indent 2\times 9^n \left(55223-94875\times 2^{1+n}+480487\times 4^n\right)\big)\,.\label{T2}
		\end{align}
	\end{small}
	Moreover, 
	\begin{align}
	&\lim_{n \to \infty}H_{\rm FO}(n) = \frac{25}{84}\,,\label{S03}\\
	&\lim_{n \to \infty}H_{\rm SO} (n)= \frac{25}{432}(N_n)^{ (\log_3 4)-1} .\label{S04}
	\end{align}
\end{theorem}
\begin{IEEEproof}
	Since $H_{\rm FO}(n)=S_n/(2N_n)$  and  $H_{\rm SO}(n)=T_n/(2N_n)$, combining  the above-obtained related quantities  and using~\eqref{S01} and~\eqref{S02}, we obtain \eqref{T1} and \eqref{T2}, which  lead to \eqref{S03} and \eqref{S04} for sufficiently large $n$.
\end{IEEEproof}

Notice that \eqref{T1} and \eqref{S03} were previously obtained in~\cite{YiZhPa20} by using a different technique.  However,  \eqref{T2} and \eqref{S04} are novel.   Theorem~\ref{TH01} shows that in large networks $G_n$, the  first-order coherence  $H_{\rm FO}(n)$ does not depend on $n$, thus not on  $N_n$,   while the second-order coherence  $H_{\rm SO}(n)$  increases sublinearly  with $N_n$.


%

\section{Network Coherence in Sierpi\'{n}ski Gaskets}

In the previous section, we  obtained an explicit formula of network coherence $H_{\rm FO}(n)$ for PSFW $G_n$, and showed that  $H_{\rm FO}(n)$ is a sublinear function of  $N_n$.  To unveil the effect of scale-free and small-world topologies on  the scaling of network coherence, in this section  we derive an analytical expression for network coherence in Sierpi\'{n}ski  gaskets with the same numbers of vertices and edges as those of PSFWs. We will show that the leading scalings of both first-order and second-order network coherence for $S_n$ are significantly different from those associated with $G_n$.



\begin{figure}
	\centering
	\includegraphics[width=0.4\textwidth]{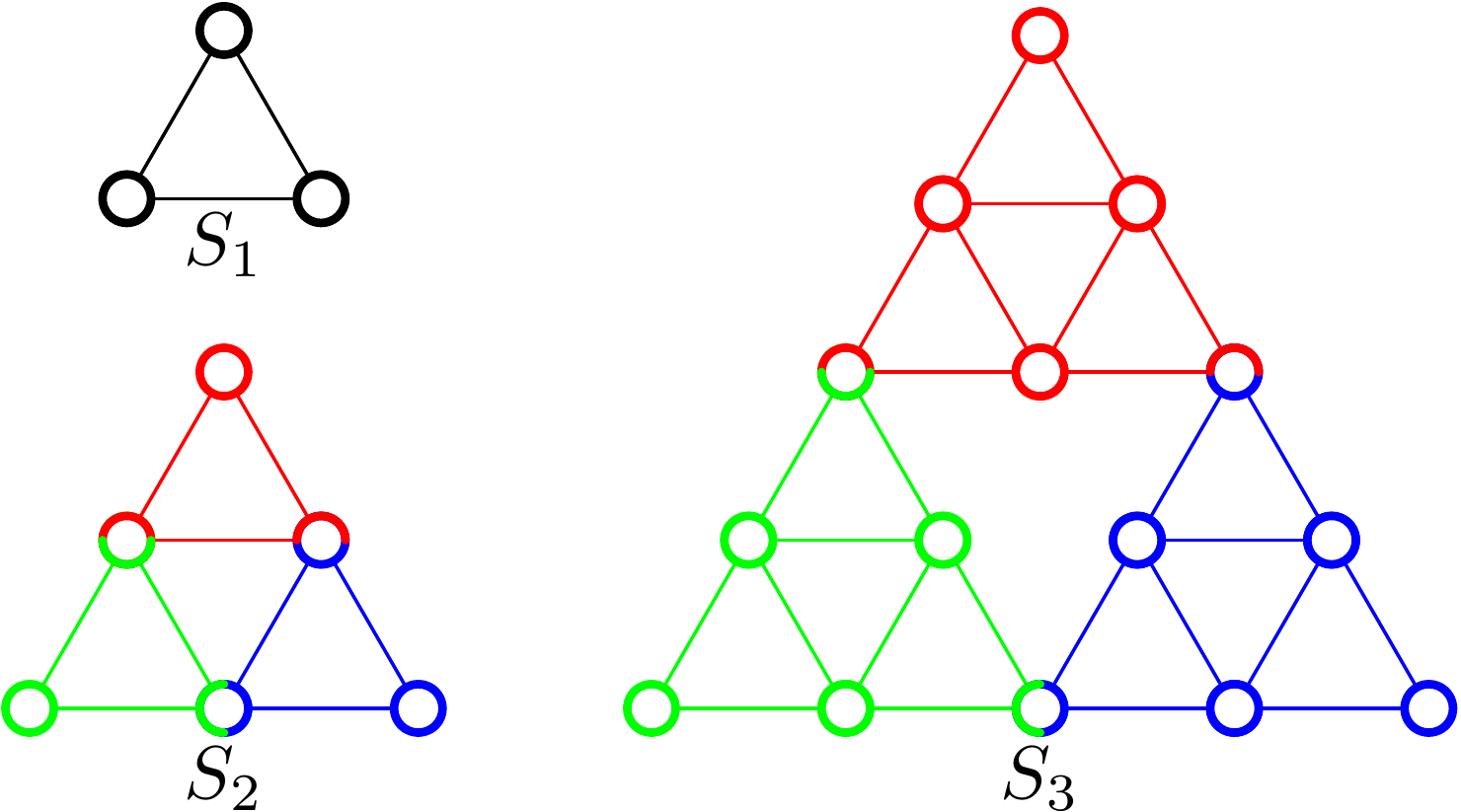}
	\caption{The first three generations of Sierpi\'{n}ski  gaskets.}
	\label{Fig.s1}
\end{figure}

The Sierpi\'{n}ski  gaskets are also iteratively constructed. Let $S_n$ ($n\ge 0$)  denote the networks after $n$ iterations. Then, the Sierpi\'{n}ski  gaskets are generated as follows.  When $n=0$,  $S_0$ is an equilateral triangle with three vertices and three edges. For $n=1$, the three edges of the equilateral triangle $S_0$ are bisected and the central triangle is removed, yielding  $S_1$ containing three copies of the original triangle. For  $n\ge 1$,  $S_n$ is generated from $S_{n-1}$ by performing the  bisection and removal  procedure for each upward pointing triangle in  $S_{n-1}$.   Fig.~\ref{Fig.s1} illustrates the first three generations of Sierpi\'{n}ski  gaskets.

Both the number of vertices and the number of edges in $S_n$ are the same as those of $G_n$. That is, there are $N_n=\frac{3^{n+1}+3}{2}$  vertices and  $E_n=3^{n+1}$ edges in $S_n$.  Many other properties of $S_n$ and $G_n$ are quite  different from each other. For example, Sierpi\'{n}ski  gaskets are neither scale-free nor small-world. They are homogeneous, with the degrees of the three outmost vertices  equal to 2, while the the degrees of other vertices being 4.

\begin{figure}
	\centering
	\includegraphics[width=0.5\textwidth]{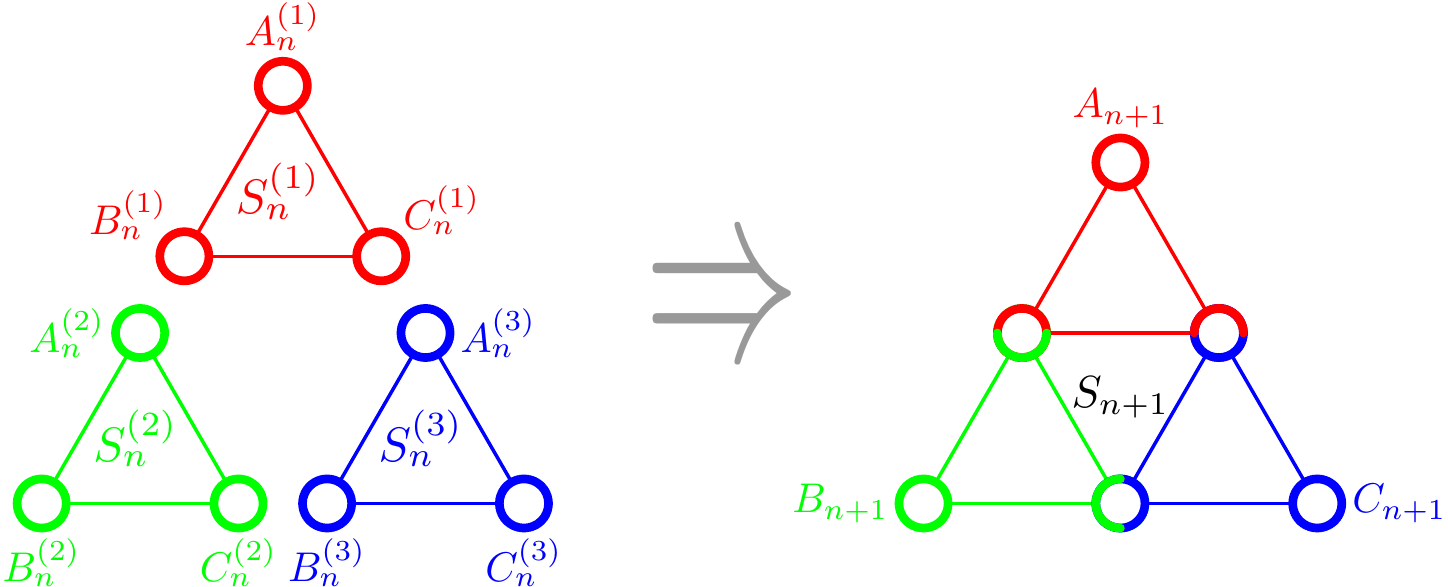}
	\caption{Second construction of  Sierpi\'{n}ski  gaskets, highlighting the self-similar structure.}
	\label{Fig.s2}
\end{figure}

Despite the difference between $S_n$ and $G_n$,  there are some similarity between them. For instance,  both graphs have cycles of various lengths. Moreover, Sierpi\'{n}ski  gaskets are likewise self-similar, as can be seen from the following alternative construction approach.  We denote  the three outmost vertices in  $S_n$ with degree 2 by $A_n$, $B_n$, and $C_n$, respectively.  Then $S_{n+1}$ is obtained from $S_n$  by joining three copies  of $S_n$ at their outmost vertices, as shown in Fig.~\ref{Fig.s2}. Let $S_n^{(\theta)}$, $\theta=1,2,3$,  represent the three replicas of $S_n$, with outmost vertices  $A_n^{(\theta)}$, $B_n^{(\theta)}$, and $C_n^{(\theta)}$. Then, $S_{n+1}$ is created  by coalescing $S_n^{(\theta)}$, $\theta=1,2,3$, with  $A_n^{(1)}$, $B_n^{(2)}$, and $C_n^{(3)}$ being the three outmost vertices of $S_{n+1}$ .

Let $H_{\rm FO}(n)$  and  $H_{\rm SO}(n)$  denote, respectively, the first-order network coherence and second-order network coherence for  $S_n$. By using a similar method and procedure, we can obtain  exact solutions for $H_{\rm FO}(n)$  and  $H_{\rm SO}(n)$ and the leading scalings for  $S_n$, as summarized in the following theorem.
\begin{theorem}\label{TH02}
	For the  Sierpi\'{n}ski   gasket $S_n$ with $n \geq 1$, the first-order coherence $H_{\rm FO}(n)$ and the second-order coherence $H_{\rm SO}(n)$ of the system with dynamics in \eqref{Hfo} and \eqref{Hso} are
	\begin{align}
	&H_{\rm FO}(n)=\frac{1}{20\times 3^{2+n} \left(1+3^n\right)^2}\big(4\times 3^n+\nonumber\\
	&\setlength{\parindent}{4.8em}\indent 2\times 3^{1+2 n}-3^{1+3 n}+13\times 3^{1+n} 5^n+\nonumber\\
	&\setlength{\parindent}{4.8em}\indent 4\times 5^{1+n}+14\times 45^n\big)\,,\label{T01}\\
	&H_{\rm SO}(n)=\frac{1}{400 \left(1+3^n\right)^3 9^{2+n}}\big(86\times 3^{1+4 n}-\nonumber\\
	&\setlength{\parindent}{4.8em}\indent 2\times 3^{2+5 n}+754\times 3^{1+2 n} 5^n+568\times \nonumber\\
	&\setlength{\parindent}{4.8em}\indent 3^{1+3 n} 5^n+32\times 3^{2+n} 5^{1+2 n}+119\times 9^n+\nonumber\\
	&\setlength{\parindent}{4.8em}\indent 28\times 5^n 9^{1+2 n}+ 64\times 15^{1+n}+8\times \nonumber\\
	&\setlength{\parindent}{4.8em}\indent 9^{1+2 n} 25^n+24\times 25^{1+n}+320\times 27^n+\nonumber\\
	&\setlength{\parindent}{4.8em}\indent 1237\times 225^n+394\times 675^n\big)\,.\label{T02}
	\end{align}
	Moreover,
	\begin{eqnarray}\label{T03}
	\lim_{n \to \infty}H_{\rm FO} (n)=\frac{7}{90} (N_n)^{(\log_3 5)-1}
	\end{eqnarray}
	and
	\begin{eqnarray}\label{T04}
	\lim_{n \to \infty}H_{\rm SO}(n) = \frac{1}{450}(N_n)^{(\log_3 {25})-1}\,.
	\end{eqnarray}
\end{theorem}

Theorem~\ref{TH02} shows that  the behaviors of both first-order coherence  $H_{\rm FO}(n)$ and second-order coherence $H_{\rm SO}(n)$ in Sierpi\'{n}ski gaskets significantly differ  from those of PSFWs.  For large Sierpi\'{n}ski   gaskets  $S_n$,  the  first-order coherence  $H_{\rm FO}(n)$ increases sublinearly  with  $N_n$,  while the second-order coherence  $H_{\rm SO}(n)$  behaves superlinearly  with $N_n$.

\begin{figure}
	\centering
	\includegraphics[width=0.52\textwidth,trim=10 10 0 40]{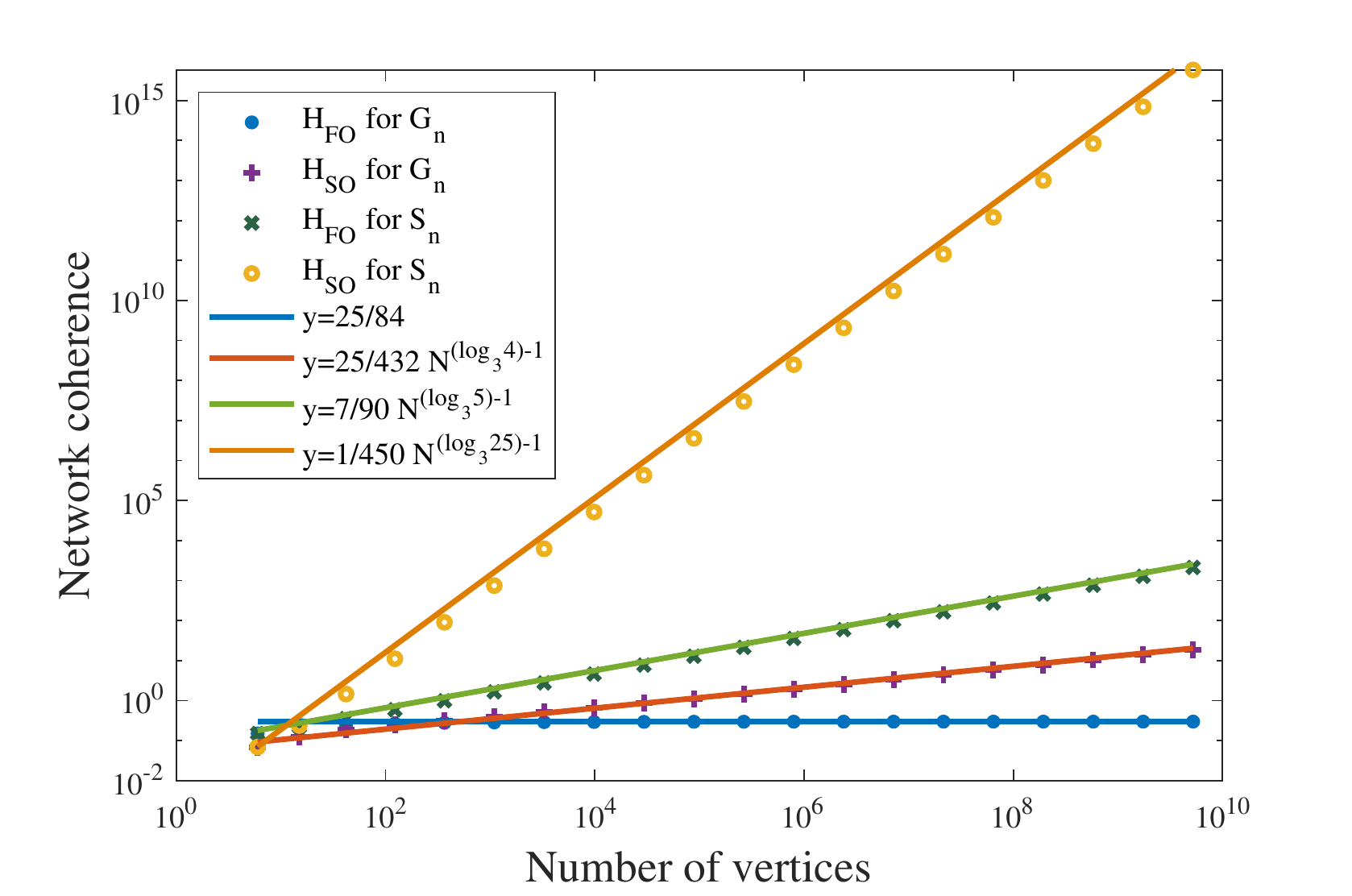}
	\caption{First-order and second-order network coherences  versus $N_n$ in both $G_n$ and $S_n$ on a log-log scale, with  $n$ changing from 1 to 20. The exact results of $G_n$ are  calculated by \eqref{T1} and \eqref{T2}, while the explicit results of $S_n$ are  obtained by \eqref{T01} and \eqref{T02}.  The approximate results are obtained by \eqref{S03} and \eqref{S04} for $G_n$, and by \eqref{T03} and \eqref{T04} for $S_n$.  }
	\label{Fig.4}
\end{figure}


In Fig.~\ref{Fig.4}, we report a direct comparison of approximate and exact results about first-order coherence $H_{\rm FO}(n)$ and second-order coherence $H_{\rm SO}(n)$ in PSFWs $G_n$ and Sierpi\'{n}ski  gaskets $S_n$ for various $n$. For moderately large $n$,   the exact  and approximate results agree with each other.

\section{Result Analysis}\label{analysis.sec}

In the previous sections, we have investigated the noisy second-order consensus dynamics of some real-life scale-free networks and a class of scale-free model networks. For all these studied scale-free networks, their second-order network coherence scales sublinearly with the number of vertices, $N$. Note that for a network, its second-order coherence is completely determined by the sum of the square reciprocal of every non-zero eigenvalue of its Laplacian matrix. Since for scale-free networks, the eigenvalues and their distributions are closely related to the network structures~\cite{ZhChYe10}, the sublinear scaling for coherence observed for the considered scale-free networks lies in their intrinsic structural characteristics, particularly the scale-free small-world topology and cycles of different lengths.

As shown in~\cite{YiYaZhPa18}, the second-order coherence of a network is determined by the average of biharmonic distances $\Theta_{ij}$ over all pairs of vertices. By formulas~\eqref{rdis} and~\eqref{bdis}, for any pair of vertices $i$ and $j$, both resistance distance $\Omega_{ij}$ and biharmonic distance $\Theta_{ij}$ are combinations of $\frac{1}{\lambda_k}(u_{ki}-u_{kj})^2$, $k=1,2,\ldots,N-1$. For the resistance distance, the weight of each term is 1, while for the biharmonic distance, the weight varies, with a  larger term $\frac{1}{\lambda_k}(u_{ki}-u_{kj})^2$ corresponding to a larger weight $\frac{1}{\lambda_k}$.  Thus, for most graphs, $\Theta_{ij}$ is greater than $\Omega_{ij}$. In a scale-free network, the existence of large-degree vertices connected to many other vertices is accompanied by the small-world property, characterized by at most a logarithmical growth average path length~\cite{Ne03}. Moreover, for a scale-free small-world network, its average resistance distance is even smaller, converging to a constant~\cite{YiZhPa20}. In contrast to the constant average resistance distance, the average of biharmonic distances is dependent on $N$, scaling sublinearly with $N$. Next, we show that this sublinear scaling is an aggregation of scale-free, small-world and loopy structures, since neither power-law small-world behavior nor cycles alone can ensure a sublinear coherence, but it leads to a linear or superlinear scaling.

It was reported that for the scale-free small-world Koch network~\cite{YiZhShCh17}, its second-order network coherence $H_{\mathrm{SO}}$ behaves linearly with the number of vertices, $N$, which is also observed for the small-world hierarchical graphs~\cite{QiZhYiLi19} with an exponential degree distribution. Both Koch networks and hierarchical graphs are highly clustered, but have only small cycles such as triangles, lacking cycles of various lengths. Thus, the existence of cycles of different lengths is necessary for sublinear scaling of $H_{\mathrm{SO}}$ in a network. However, cycles do not suffice to guarantee a sublinear scaling of $H_{\mathrm{SO}}$. For example, in Sierpin\'ski gaskets, there are cycles of various lengths, but their $H_{\mathrm{SO}}$ is a superlinear function of $N$ as given by formula~\eqref{T04}. This, in turn, indicates that scale-free small-world topology is only necessary for a sublinear scaling of $H_{\mathrm{SO}}$ in a sparse network.

Note that although the small-world is an accompanying phenomenon of the power-law behavior~\cite{Ne03}, small-world and loopy structures cannot lead to the sublinear scaling of second-order network coherence. For example, using the result in~\cite{YiZhLiCh15} and the technique in this current paper, we can determine the analytical expression of second-order coherence for the Farey graphs, which scales linearly with the number of nodes, being quite different from the sublinear scaling observed for the PSFWs. By construction, Farey graphs are subgraphs of PSFWs, both of which are small-world and highly clustered, with many cycles at different scales. The reason for the scaling distinction of the second-order coherence between Farey graphs and PSFWs lies in, at least partially, the scale-free topology of PSFWs that the Farey graphs do not possess.

\section{Conclusion}\label{conclusion.sec}

A large variety of real-world networks are sparse and loopy, and exhibit simultaneously scale-free and small-world features. These structural properties have a substantial influence on different dynamics running on such networks. In this paper, we presented an extensive study on second-order consensus in noisy networks with these properties, focusing on its robustness measured by network coherence that is characterized by the average steady-state variance of the system. We first studied numerically the network coherence for some representative real scale-free networks, which grows sublinearly with the vertex number $N$. We then determined exactly the coherence for a class of deterministic scale-free networks, PSFWs, which is also a sublinear function of $N$. Moreover, we studied analytically the coherence for  Sierpin\'ski gaskets with the same numbers of vertices and edges as the PSFWs, the leading scaling of which scales superlinearly in $N$.  We concluded that the scale-free, small-world, and loopy structures are responsible for the observed sublinear scaling of coherence for the studied networks.

It should be mentioned that we only addressed second-order noisy consensus on undirected binary networks, concentrating on scale-free, small-world, and loopy properties on the effects of network coherence. Future work should include the following directions. First, it would be of interest to consider second-order noisy consensus on directed~\cite{LiLuHu18,LuLiLi20} and weighted~\cite{BaBaPaVe04, QiLiZh18} communication graphs, with an aim to explore the influences of one-way action or distribution of edge weights on network coherence. Another direction is to examine second-order noisy linear consensus networks in the presense of time-delay~\cite{SoGhMo19,LuliXi20}. Moreover, of particular interest is to consider the case that both scalar-valued states of each agent are subject to disturbances. Finally, our method and process for computing the network coherence are only applicable to deterministically growing self-slimilar networks, it is of great significance to modify or extend them to more general networks.

\appendices

\input{2ndConsensus.bbl}

\bibliographystyle{IEEEtran}
\begin{IEEEbiography}[{\includegraphics[width=1in,height=1.25in,clip,keepaspectratio]{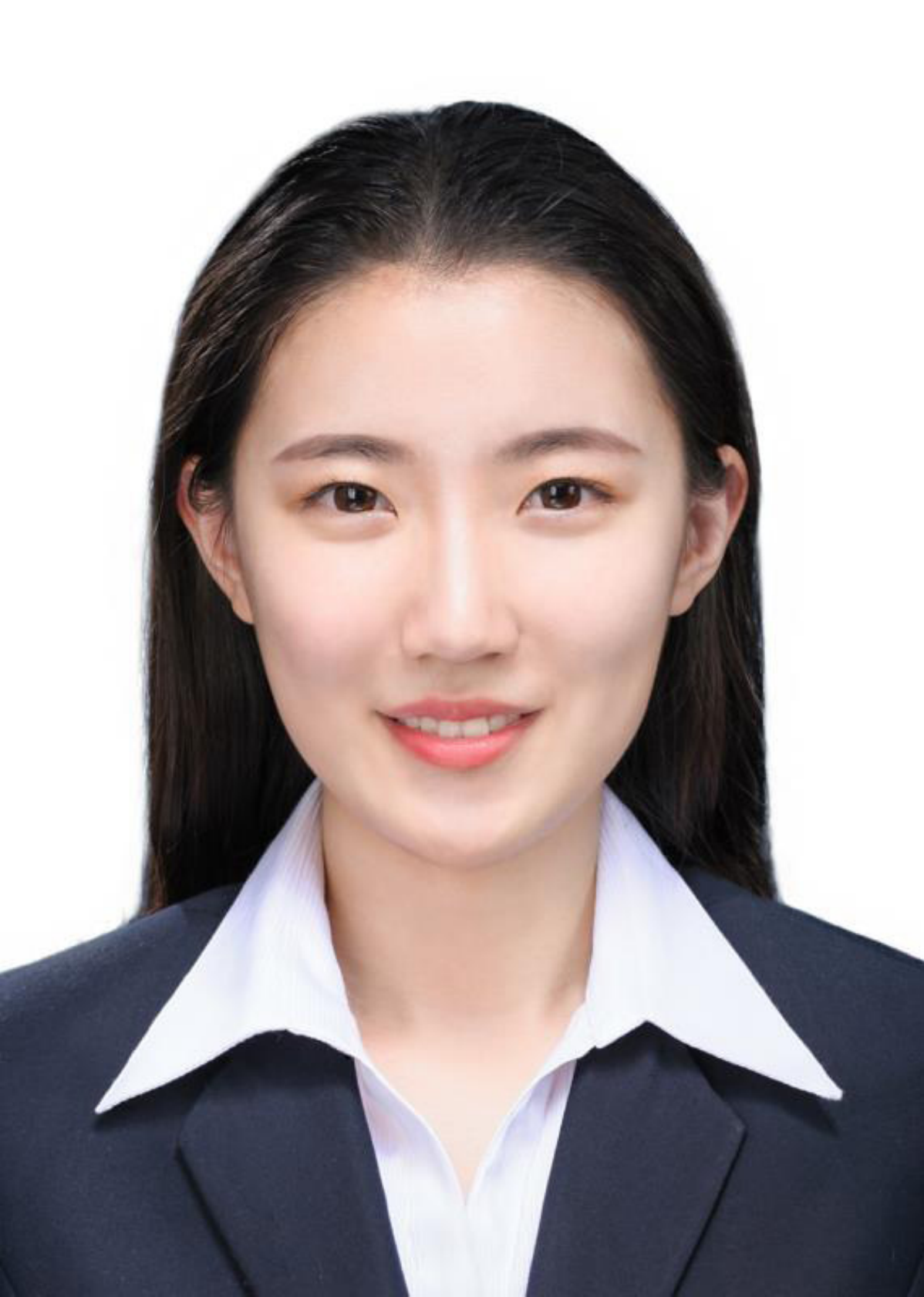}}]{Wanyue Xu}
	received the B.Eng. degree in computer science and technology, Shangdong University, Weihai, China, in 2019. She is currently pursuing the Master degree in the School of Computer Science, Fudan University, Shanghai, China. Her research interests include network science, graph data mining, social network analysis, and random walks.
	Ms. Xu has published several papers in international journals or conferences, including TCYB, WWW, WSDM, and ICDM.
\end{IEEEbiography}
\begin{IEEEbiography}[{\includegraphics[width=1in,height=1.25in,clip,keepaspectratio]{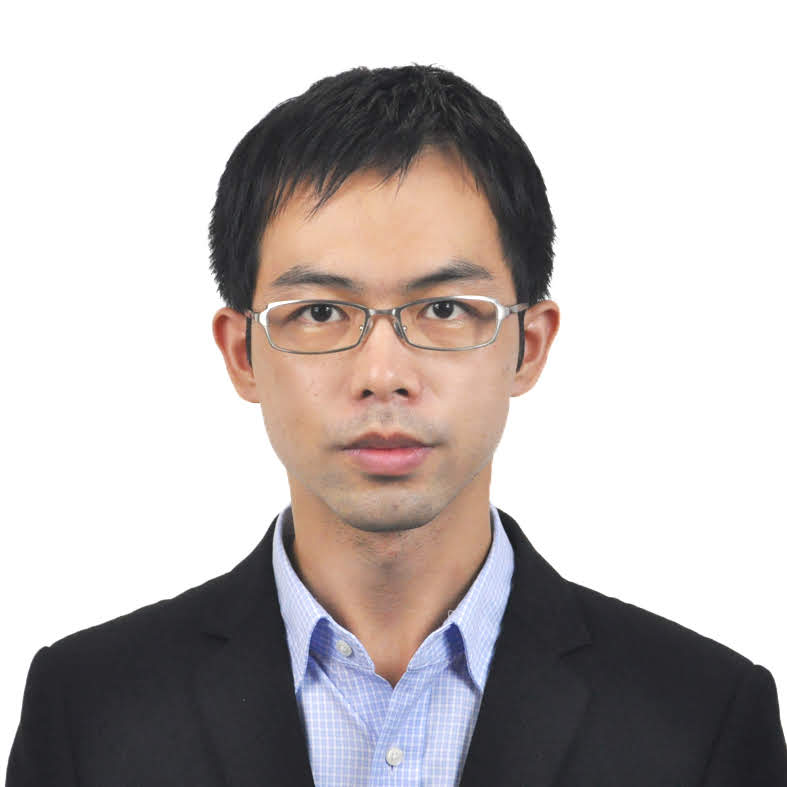}}]{Bin Wu}
	received the B.S. degree and the M.Sc. degree in computer science from Fudan University, Shanghai, China, in 2011 and 2014, respectively. His research interests include complex networks, random walks, and spectral graph theory.
\end{IEEEbiography}
\begin{IEEEbiography}[{\includegraphics[width=1in,height=1.25in,clip,keepaspectratio]{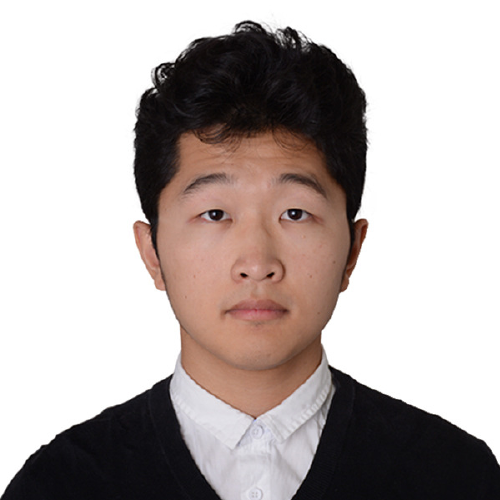}}]{Zuobai Zhang}
	 is currently an undergraduate student working toward the B.S. degree in the School of Computer Science, Fudan University, Shanghai, China. His research interests include graph algorithms, social networks, and network science.\\
	Mr. Zhang has published several papers in international journals or conferences, including WSDM, WWW and IEEE Transactions on Cybernetics. He won a Silver Medal in National Olympiad in Informatics of China in 2016 and several Gold Medals in ICPC Asia Regional Contests.
\end{IEEEbiography}
\begin{IEEEbiography}[{\includegraphics[width=1in,height=1.25in,clip,keepaspectratio]{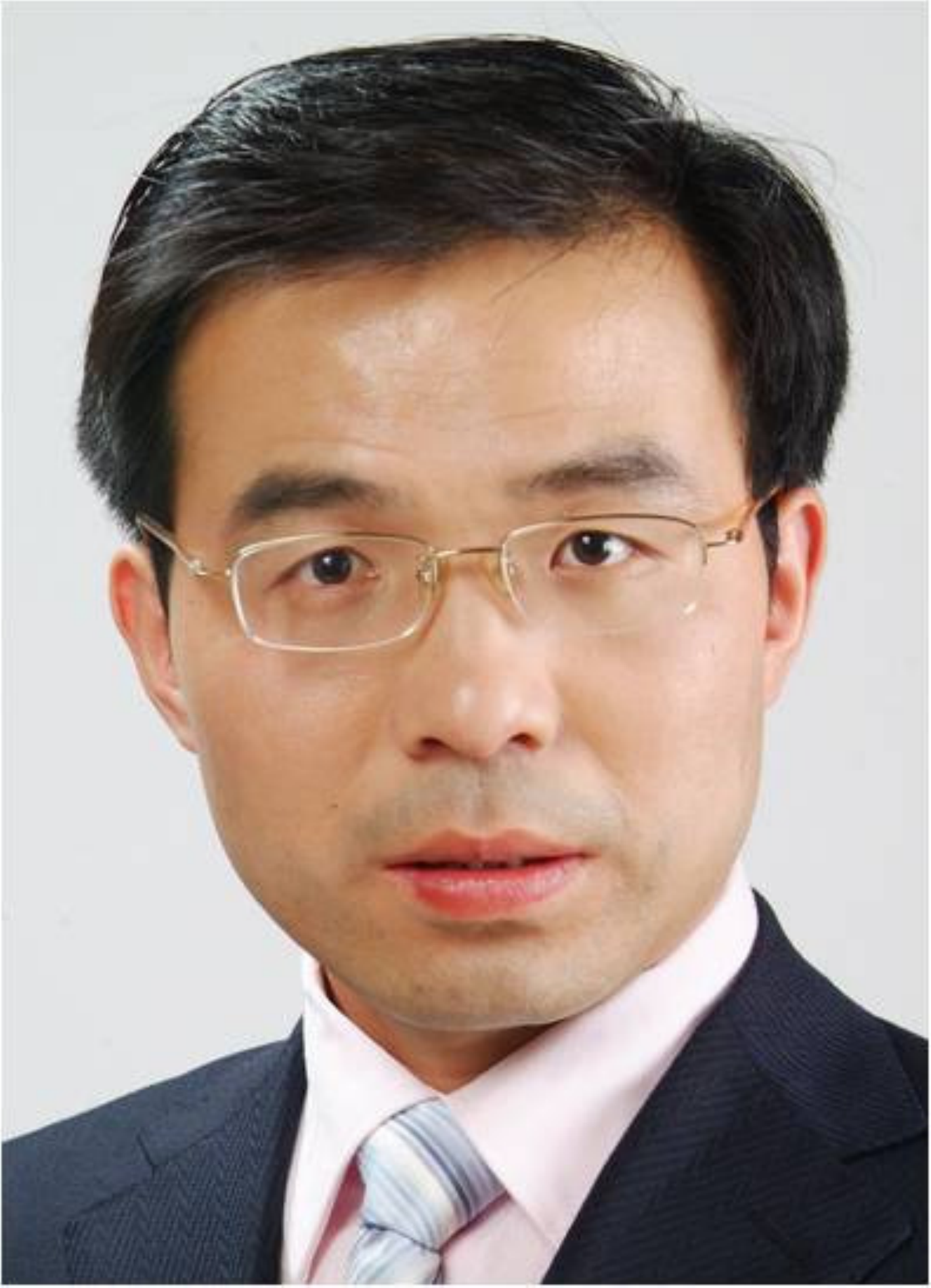}}]{Zhongzhi Zhang}
	(M'19) received the B.Sc. degree in applied mathematics from Anhui University, Hefei, China, in 1997 and the Ph.D. degree in management science and engineering from  Dalian University of Technology, Dalian, China, in 2006. \\
	From 2006 to 2008, he was a Post-Doctoral Research Fellow with Fudan University, Shanghai, China, where he is currently a Full Professor with the School of Computer Science. He has published over 140 papers in international journals or conferences. He has over 3000 ISI Web of Science citations with an H-index of 32 according to the Clarivate. He was one of the most cited Chinese researchers (Elsevier) in 2019. His current research interests include network science, graph data mining, social network analysis, spectral graph theory, and random walks. \\
	Dr. Zhang was a recipient of the Excellent Doctoral Dissertation Award of Liaoning Province, China, in 2007, the Excellent Post-Doctor Award of Fudan University in 2008, the Shanghai Natural Science Award (third class) in 2013, and the Wilkes Award for the best paper published in The Computer Journal in 2019. He is a member of the IEEE.
\end{IEEEbiography}
\begin{IEEEbiography}[{\includegraphics[width=1.1in,height=2in,clip,keepaspectratio]{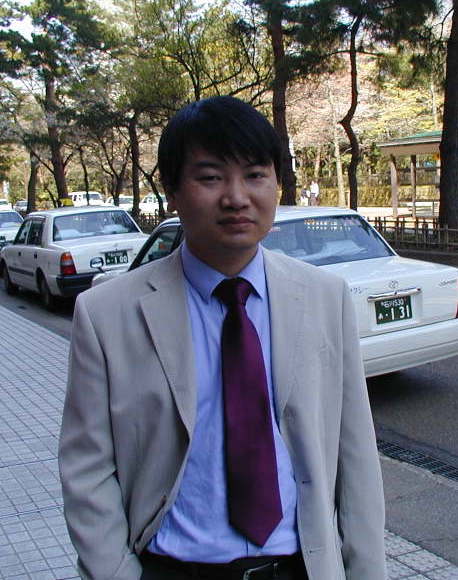}}]{Haibin Kan}
	received the Ph.D. degree from Fudan University, Shanghai, China, in 1999. Then, he became a faculty of Fudan University. From June 2002 to February 2006, he was with the Japan Advanced Institute of Science and Technology as an assistant professor. He returned to Fudan University in February 2006, where he is currently a full professor. His research interests include coding theory, cryptography, and computation complexity.
\end{IEEEbiography}
\begin{IEEEbiography}[{\includegraphics[width=1in,height=1.25in,clip,keepaspectratio]{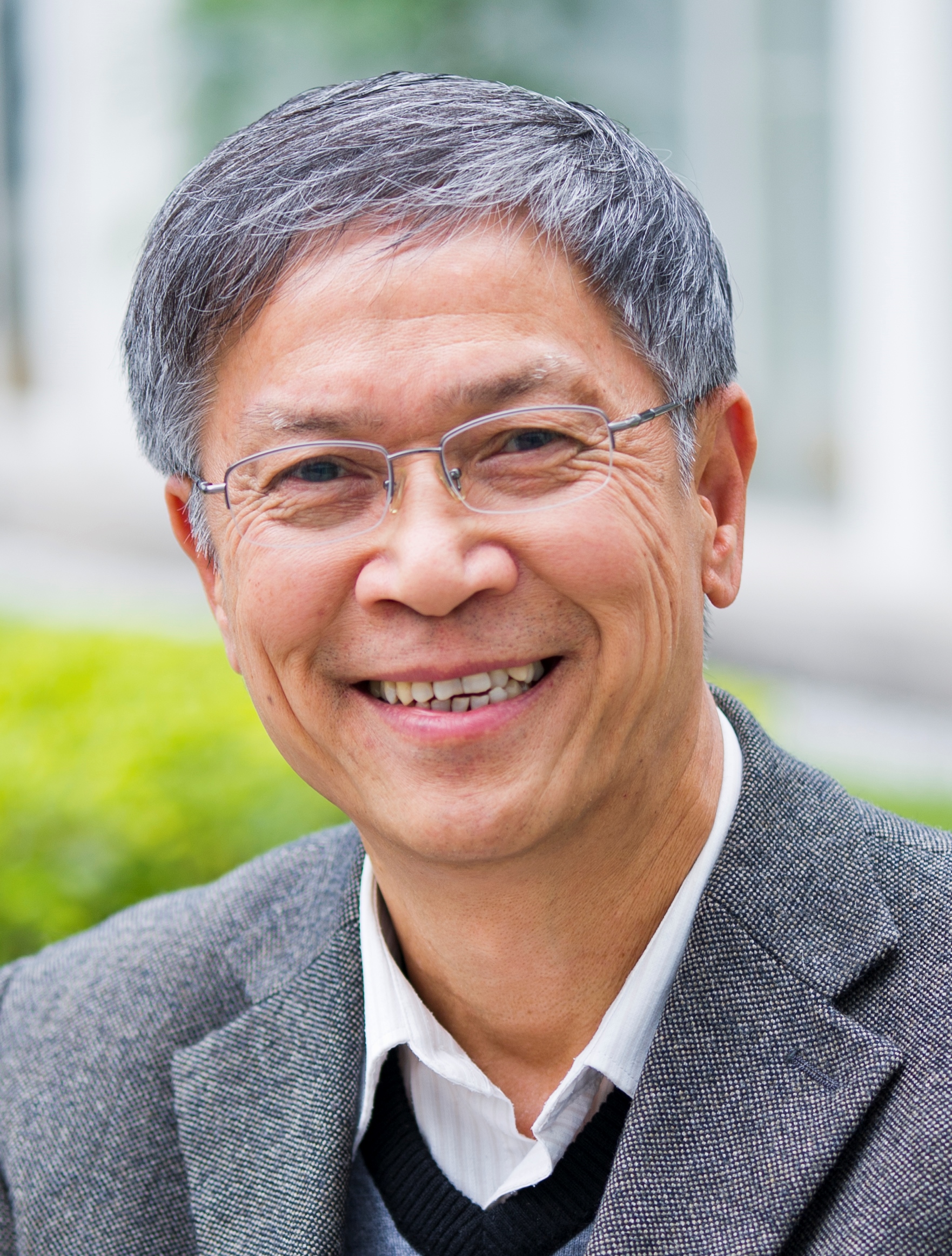}}]{Guanrong Chen}
	 (M’89–SM’92–F’97) received the M.Sc. degree in computer science from Sun Yat-sen University, Guangzhou, China, in 1981, and the Ph.D. degree in applied mathematics from Texas A\&M University, College Station, TX, USA, in 1987. \\
	 He was a Tenured Full Professor with the University of Houston, Houston, TX, USA. He has been a Chair Professor and the Founding Director of the Center for Chaos and Complex Networks, City University of Hong Kong, Hong Kong, since 2000. \\
	 Dr. Chen is a member of the Academia Europaea and a fellow of The World Academy of Sciences. He was a recipient of the 2011 Euler Gold Medal from Russia and conferred Honorary Doctorate by Saint Petersburg State University, Russia, in 2011, and by the University of Le Havre, Nomandie, France, in 2014. He is a Highly Cited Researcher in engineering and also in mathematics according to Thomson Reuters.
\end{IEEEbiography}
\end{document}

%% file: 2ndConsensus.bbl
\providecommand{\noopsort}[1]{}\providecommand{\singleletter}[1]{#1}